\documentclass[12pt]{scrartcl}

\usepackage{amsthm, amsmath, amssymb}
\usepackage[alphabetic]{amsrefs}
\usepackage{graphicx}
\usepackage{float}
\usepackage{enumitem}
\allowdisplaybreaks

\newtheorem{theorem}{Theorem}

\newtheorem{proposition}[theorem]{Proposition}

\theoremstyle{definition}

\newcommand{\R}{\mathbb{R}}

\newcommand{\Per}{\mbox{Per}}

\renewcommand{\epsilon}{\varepsilon}

\renewcommand{\leq}{\leqslant}
\renewcommand{\le}{\leqslant}
\renewcommand{\geq}{\geqslant}

\makeatletter

\title{Enhanced boundary regularity \\ of planar nonlocal minimal graphs, \\
and a butterfly effect\thanks{This work has been supported by the Australian
Research Council Discovery Project DP170104880
``N.E.W. Nonlocal Equations at Work'',
by the
DECRA Project DE180100957 ``PDEs, free boundaries
and applications'', by
the European Research Council Project
``G.H.A.I.A. Geometric and Harmonic Analysis with Interdisciplinary Applications''
and by
the INdAM-DP-COFUND-2015
Grant 713485 ``Doctoral Programme in
Mathematics and/or Applications Cofunded by Marie
Sklodowska-Curie Actions''. 
Part of this work was written on the occasion of a pleasant and fruitful visit
of the second and third authors to the Department of Mathematics
and Statistics of the University of Western Australia, which we thank for the very
warm hospitality.}}

\author{Serena Dipierro,
Aleksandr Dzhugan,\\
Nicol\`o Forcillo
and Enrico Valdinoci}

\date{} 

\begin{document}

\maketitle 

\begin{abstract} In this note, we showcase some recent results
obtained in~\cite{2019-1} concerning the stickiness
properties of nonlocal minimal graphs in the plane.
To start with, the nonlocal minimal graphs in the plane
enjoy an enhanced boundary regularity, since
boundary continuity with respect to the external datum
is sufficient to ensure differentiability across the boundary of the domain.

As a matter of fact, the H\"older exponent of the derivative
is in this situation sufficiently high to provide the validity
of the Euler-Lagrange equation at boundary points as well.

{F}rom this, using a sliding method, one also
deduces that
the stickiness phenomenon is generic for
nonlocal minimal graphs in the plane,
since an arbitrarily small perturbation of continuous 
nonlocal minimal graphs can produce boundary discontinuities
(making the continuous case somehow ``exceptional''
in this framework).
\end{abstract}

\section*{Three questions on planar nonlocal minimal graphs}

Nonlocal minimal surfaces are a beautiful -- and
extremely challenging -- topic of research.
The novelty of the subject, together with its intrinsic cross-disciplinary
nature, requires the combination of techniques from different fields,
including calculus of variations, geometric measure theory, geometric
analysis, differential geometry, partial differential and integro-differential equations.
The solution of the problems posed by this intriguing scenario
is usually based on brand new approaches and opens several perspectives
in both pure and applied mathematics.
\medskip

Moreover, nonlocal minimal surfaces offer a number of important,
and very often surprising, differences with respect to the classical case.
Among these differences,
we believe that the ones related to new
``boundary effects'' are of particular importance,
also in view of some ``stickiness phenomena'' that have been
recently discovered and which seem to play a crucial role in the understanding
of phenomena relying on long-range interactions.
The goal of this note is to recall some recent results in this direction,
and to describe the peculiar boundary situation exhibited by
planar nonlocal minimal graphs.\medskip

To this end, we recall
the definition of $s$-perimeter introduced in \cite{MR2675483}. Namely,
given~$s\in (0,1)$ and two measurable, disjoint
sets~$A$, $B\subseteq \R^n$,
we define the nonlocal set-interaction as
\[I(A,B):=\iint\limits_{A\times B}\frac{dx\hspace{0.05cm}dy}{\left|x-y\right|^{n+s}}.\]
Also, if~$\Omega\subset\R^n$ is a bounded set with Lipschitz
boundary, and $E\subseteq\R^n$ is a measurable set,
we define the $s$-perimeter of $E$ in $\Omega$ as
\[ \Per_s(E,\Omega)=I(E\cap\Omega,E^c\cap\Omega)+I(E\cap\Omega,E^c\cap \Omega^c)+I(E\cap \Omega^c,E^c\cap \Omega).\]
The name of $s$-perimeter for this type of functionals
is motivated by the fact that, as~$s\nearrow1$, this
functional recovers the classical notion of perimeter
(in various forms,
including functional estimates, $\Gamma$-convergence,
density estimates, clean ball conditions, isoperimetric inequalities, etc.,
see~\cite{MR1942116, MR1945278, MR1942130, MR2765717, MR2782803}).
On the other hand, as~$s\searrow0$, the $s$-perimeter
is related to suitable weighted Lebesgue measures,
in which the weights take into account the behavior of the set at infinity
(see~\cite{MR1940355, MR3007726}), and these features already somewhat
suggest that the problem for~$s$ close to~$1$ may
be more ``regular'' and ``close to classical variational problems''
than the problem for~$s$ close to~$0$.\medskip

We say that~$E$ is $s$-minimal in~$\Omega$ if
$$ \Per_s(E,\Omega)\leq\Per_s(F,\Omega)$$
for every $F\subset\R^n$ such that $F\setminus\Omega=E\setminus\Omega$.

If~$\tilde\Omega\subseteq\R^n$ is unbounded, one can also say that~$E$
is $s$-minimal in $\tilde\Omega$ if $E$ is $s$-minimal
in~$\Omega$,
for all bounded Lipschitz sets $\Omega\Subset\tilde\Omega$.
We refer to~\cite{2018} for a comprehensive description
of these minimization problems.\medskip

The interior regularity theory of $s$-minimizers
is an important topic of contemporary investigation,
and complete results are available only in the plane, or when the fractional
parameter~$s$ is sufficiently close to~$1$,
see~\cite{MR3090533, MR3107529, MR3331523}.
See also~\cite{MR3981295} for quantitative bounds and
regularity results
of $BV$ type for stable solutions.\medskip

The theory of nonlocal minimal surfaces is also related
to nonlocal isoperimetric problems (see~\cite{MR2425175, MR2469027, MR2799577, MR3322379, MR3412379, MR3640534, MR3732175}),
to fractional mean curvature equations (see~\cite{MR2487027, MR2564467, MR3485130, MR3744919, MR3770173, MR3836150, MR3881478}) and to nonlocal geometric flows
(see~\cite{MR3023439, MR3156889, MR3401008, MR3713894, MR3732178, MR3778164, MR3951024, MR4000255, Short}).
\medskip

Among all the possible minimization frameworks, the one of the graphs
seems to play a special role, since it enjoys a number of
structural features and can provide a solid guideline for the general
theory. To introduce this setting, given a measurable function~$u:\R^{n-1}\to\R$,
we use the notation
\begin{equation}\label{definition-E-u}
E_u:=\big\{ (x_1,\dots,x_{n-1},x_n)\in\R^n {\mbox{ s.t. }} x_n<u(x_1,\dots,x_{n-1})\big\}.
\end{equation}
Then, given a bounded Lipschitz domain~$\Omega_o\subset\R^{n-1}$,
we say that~$u$ is $s$-minimal in~$\Omega_o$ if~$E_u$ is $s$-minimal
in~$\Omega_o\times\R$.

The graphical case constitutes a useful building block for the general theory
since it provides a ``stable'' framework to work with, in the sense that if~$E$
is a graph outside~$\Omega_o\times\R$, then the $s$-minimizer
in~$\Omega_o\times\R$ is a graph as well, see~\cite{MR3516886}.

Also, the graphical structure poses some natural problems of Bernstein type
(see~\cite{MR3680376, PISA}) and enjoys several
special regularity features (see~\cite{MR3934589}). See also~\cite{MATLUC}
for additional properties of nonlocal minimal surfaces and nonlocal minimal graphs.\medskip

In this note, for the sake of concreteness,
we will focus on the planar\footnote{The higher
dimensional situation is structurally more complicated.
The first attempt to describe the boundary behavior of
higher dimensional nonlocal minimal surfaces can be found in~\cite{2019-2}.} case,
with the aim of highlighting the main features of
$s$-minimal graphs in a slab.
In this setting, given~$u_0:\R\to\R$,
the typical problem is to understand the geometric
properties of the minimizer~$u$ in~$(0,1)$
with~$u=u_0$ in~$\R\setminus(0,1)$.

When~$u_0:=0$, the minimizer~$u$ vanishes
identically, as proved in~\cite{MR2675483} using a maximum
principle argument (see also~\cite{CAB-C, PAG-C}
for recently introduced calibration methods).

\begin{figure}[H]
	\centering
	\includegraphics[width=9cm]{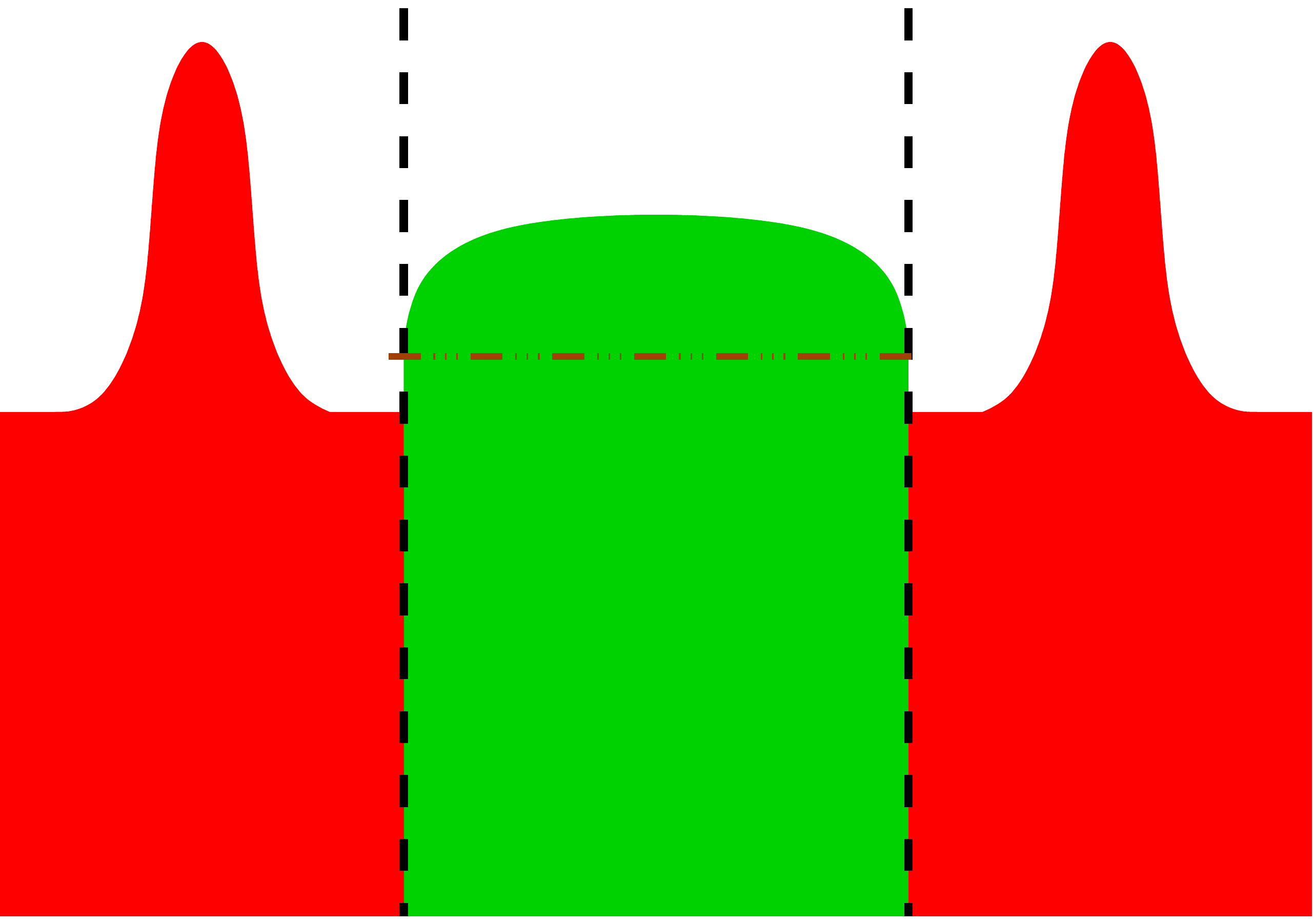}
	\caption{\textit{{Example of stickiness: initial problem with a datum with two small bumps.}}}
	\label{FIG3}
\end{figure}

An interesting question to address is what happens for
small perturbations of the exterior datum $u_0$.
The case of two small bumps was investigated
in~\cite{MR3596708}, where it was established that~$u$
remains bounded away from zero in~$(0,1)$, see Figure~\ref{FIG3}.
In particular, nonlocal minimal graphs do not necessarily meet continuously
the boundary datum -- sometimes they do, as it happens for the case~$u_0:=0$,
but small perturbations of such a datum are sufficient to produce boundary discontinuities.
Hence, the minimizing problem for nonlocal minimal graphs
is well posed in the class of functions, but not in the class of continuous functions,
since the $s$-minimal graph can turn out to be discontinuous
at the boundary (and, as a matter of fact,
this discontinuity is a jump, since the nonlocal minimal graphs
are uniformly continuous inside the domain, see~\cite{MR3516886}).\medskip

This feature is a special case of a general phenomenon that
was named ``stickiness'' in~\cite{MR3596708},
emphasizing that, differently from the classical case, nonlocal minimal surfaces
have the tendency to adhere at the domain (this may be also related
to a capillarity effect, see also~\cite{MR3707346, MR3717439} for a
specific analysis of a nonlocal capillarity theory,
and~\cite{MR3926519, CLALUC} for
several examples of sticky behaviors of $s$-minimal
surfaces).\medskip

The stickiness phenomenon detected in~\cite{MR3596708},
rather than constituting a final goal for the theory of nonlocal minimal
surfaces, served as a key to disclose a number of new directions
of investigation, including:
\begin{eqnarray*}
{\mbox{(Q1)}}&&{\mbox{\em How regular are the nonlocal minimal graphs
``coming from inside the domain''?}}\\
{\mbox{(Q2)}}&&{\mbox{\em Is the Euler-Lagrange equation coming from the variation
of the $s$-perimeter}}\\&&{\mbox{\em satisfied ``up to the boundary''?}} \\
{\mbox{(Q3)}}&&{\mbox{\em How ``typical'' is
the stickiness phenomenon?}}
\end{eqnarray*}
We will show in this note that these questions are intimately correlated
and the understanding of each of these problems sheds some light on the others.

The first results addressing~(Q1) and~(Q2)
have been obtained in~\cite{MR3532394},
in which it is shown that nonlocal minimal graphs,
in the vicinity of discontinuity boundary points, can be written
as differentiable graphs with respect to the vertical variable.
Namely, if~$u$ is $s$-minimal in~$(0,1)$
with respect to a smooth datum~$u_0:\R\to\R$
and
\begin{equation}\label{US3dj}
u(0):=\lim_{x_1\searrow0} u(x_1)> u_0(0),\end{equation}
then there exist~$\rho>0$ and~$v\in C^{1,\frac{1+s}{2}}(\R)$
such that
\begin{equation}\label{sh72} \big\{ x_1\in(0,\rho),\;x_2=u(x_1)\big\}
=\big\{ x_2\in(u(0),u(\rho)),\;x_1=v(x_2)\big\},\end{equation}
with
\begin{equation}\label{sh73}
v'( u(0))=0,
\end{equation}
and a similar statement holds true when~\eqref{US3dj} is replaced by
\begin{equation}\label{sh74}
u(0)< u_0(0).\end{equation}
We remark that, in particular, \eqref{sh72} says that~$u$
is invertible near the boundary discontinuity, and, in view of~\eqref{sh73},
\begin{equation}\label{sh744} \lim_{x_1\searrow0} u'(x_1)=+\infty.\end{equation}
With respect to question~(Q1), this says that
at boundary discontinuities the derivative of~$u$ blows up,
but the graph can be seen as the inverse of a $C^{1,\frac{1+s}{2}}$-function~$v$
which has a critical point in correspondence to the jump of~$u$.
\medskip

This fact can be used to provide a first answer to~(Q2)
at boundary discontinuities, since one can equivalently write the Euler-Lagrange
equation ``along the graph of~$v$'', and then pass it to
the limit using the regularity of~$v$
(roughly speaking, the Euler-Lagrange equation involves
a fractional curvature which is an object of order~$1+s$; then,
since~$1+\frac{1+s}2>1+s$, a control in~$C^{1,\frac{1+s}{2}}$
is sufficient to pass the equation to the limit). In this way,
one obtains that the Euler-Lagrange equation
is satisfied along the closed curve
\begin{equation}\label{CURVA}
{\mathcal{C}}:=\overline{(\partial E_u)\cap ((0,1)\times\R)}\end{equation}
provided that the solution has jump discontinuities at~$x_1=0,1$,
see Theorem B.9
in~\cite{MR3926519} for a precise statement.\medskip

After these preliminary considerations, it remains to
address~(Q1) and~(Q2) at points of boundary continuity
(this, as we will see, will also provide an answer to~(Q3)).
\medskip

The main result for~(Q1) is that
a continuous $s$-minimal graph is necessarily differentiable across the boundary
(and, in fact, of class~$C^{1,\frac{1+s}2}$).
Indeed, as proved in~\cite{2019-1}, we have that:

\begin{theorem}[Enhanced boundary regularity for planar nonlocal
minimal graphs:
continuity implies differentiability]\label{FUN}
	Let $\beta\in (s, 1)$ and ~$u_0 : \R \to \R$,
with \begin{equation} \label{u0reg}u_0 \in C^{1,\beta}([-h, 0])\end{equation}
	for some $h >0$.
	Assume that $u$ is $s$-minimal in $(0,1)$ with datum~$u_0$, and that
	\begin{equation}\label{continuity-condition}
	u(0):=\lim_{x_1\searrow0}u(x_1) = u_0(0 ).
	\end{equation}
	Then, $u \in C^{1,\gamma}([-h, 1/2])$, with
	\begin{equation}\label{definition-of-gamma}
	\gamma := \min \left\{\beta, \frac{1 + s}{2}\right\}.
	\end{equation}
\end{theorem}

\begin{figure}[H]
		\centering
		\includegraphics[width=8cm]{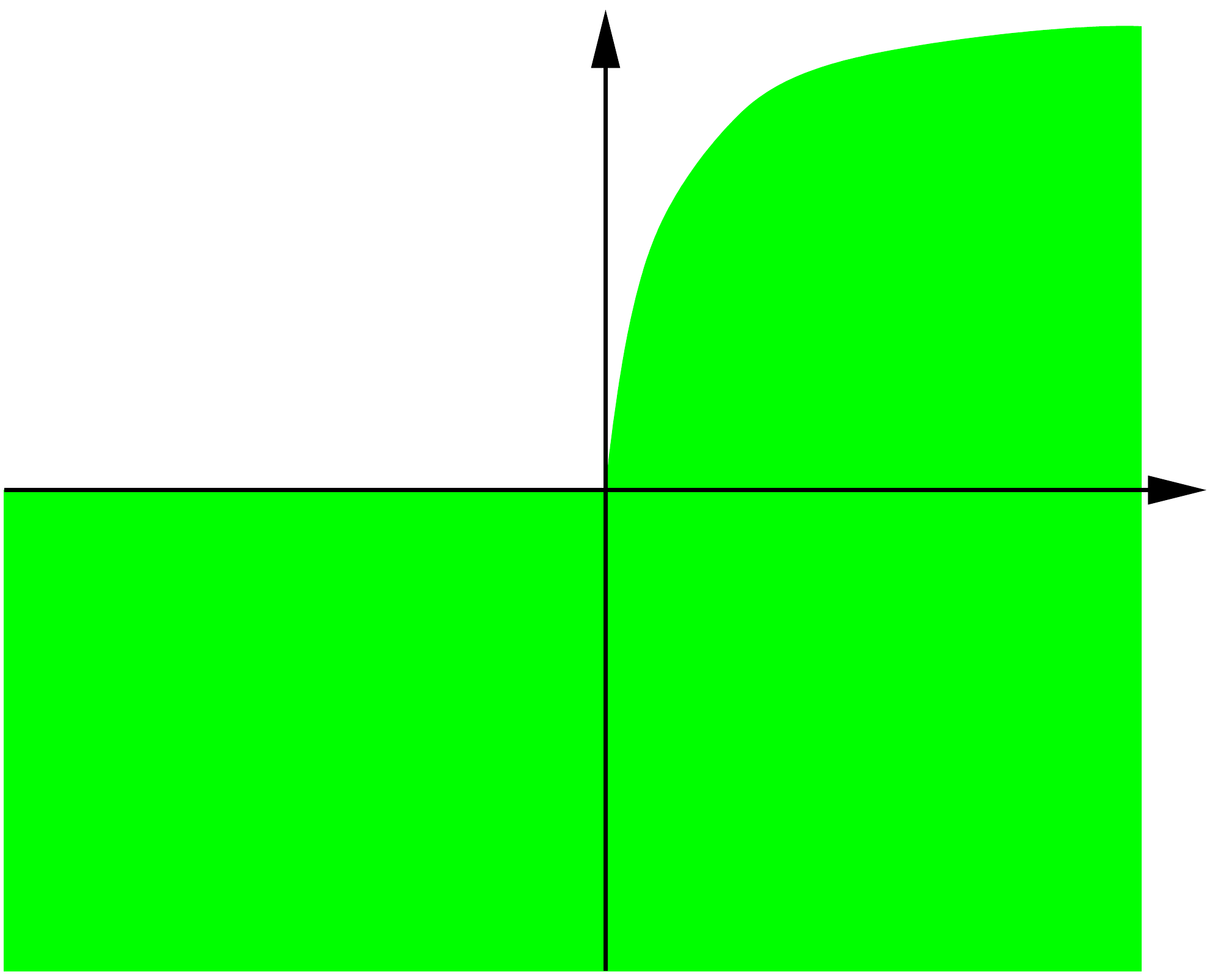}
		\caption{\textit{{Boundary behavior of an $s$-harmonic function in $(0,1)$
					vanishing in $(-1,0]$, e.g.~$u(x_1)=(x_1)_+^s$.}}}
		\label{FIG5}
	\end{figure}

When compared to
the theory of fractional linear equations,
the result in Theorem~\ref{FUN} is quite surprising since it says that continuity
is sufficient for differentiability. This is in sharp contrast with the regularity
of solutions of fractional Laplace equations such as
$$ \begin{cases}
(-\Delta)^s u =f & {\mbox{ in }}\Omega,\\
u=0&{\mbox{ in }}\R^n\setminus\Omega,
\end{cases}$$
which are in general not better than H\"older continuous at the boundary,
even when~$f$ is as smooth as we wish (see Figure~\ref{FIG5},
as well as~\cite{MR3168912} for a thorough discussion
of the boundary regularity).

For our purposes, it is interesting to observe that~\eqref{US3dj}, \eqref{sh74} and~\eqref{continuity-condition}
exhaust all the possible boundary behaviors, and
the results in~\eqref{sh72} and Theorem~\ref{FUN}
always provide a regularity of~$C^{1,\frac{1 + s}{2}}$-type up to the boundary
``in a geometric sense'': namely, 
planar
$s$-minimal graphs corresponding to {\em smooth external data}
are {\em always $C^{1,\frac{1 + s}{2}}$-curves} in the domain,
{\em up to the boundary} of the domain,
in the sense expressed by the following dichotomy:
\begin{itemize}
\item if a boundary discontinuity occurs, then the curve develops a vertical
tangent at the boundary, as given in~\eqref{sh744},
\item if the nonlocal minimal graph happens to be continuous at the boundary,
then it is actually $C^{1,\frac{1 + s}{2}}$ across the boundary.
\end{itemize}

More explicitly, we have the following result:

\begin{theorem}[Regularity of $s$-minimal curves]\label{QQ1}
Let~$u_0 : \R \to \R$, with $ u_0 \in C^{1,\frac{1+s}2}([-h, 0])$
for some $h >0$.
Assume that $u$ is $s$-minimal in $(0,1)$ with datum~$u_0$.

Then, the set ${\mathcal{C}}$ in~\eqref{CURVA}
is a $C^{1,\frac{1+s}2}$-curve.
\end{theorem}

We observe that not only Theorem~\ref{QQ1} provides a complete answer
to question~(Q1), but it also answers question~(Q2), since
one can write the Euler-Lagrange equation at the points in the interior
of the domain and then use the regularity of the curve in Theorem~\ref{QQ1}
in order to reach the boundary of the domain. In this way, one obtains that:

\begin{theorem}[Pointwise validity of the Euler-Lagrange equation]\label{VALID}
Let $\beta\in (s, 1)$ and~$u_0 : \R \to \R$, with $ u_0 \in C^{1,\beta}([-h, 0])$
for some $h >0$.
Assume that $u$ is $s$-minimal in $(0,1)$ with datum~$u_0$,
and let~${\mathcal{C}}$ be as in~\eqref{CURVA}.

Then
\begin{equation}\label{EELLQQ}
\int\limits_{\R^2} \frac{\chi_{\R^2\setminus E_u}(y)-\chi_{E_u}(y)}{|x-y|^{2+s}}\,d
y=0\end{equation}
for all $x=(x_1,x_2)\in{\mathcal{C}}$.
\end{theorem}

As customary, equation~\eqref{EELLQQ} can be considered the Euler-Lagrange
equation associated with the nonlocal perimeter and its left hand side
can be regarded as a nonlocal mean curvature (see e.g.~\cite{MR3230079, MR3733825, MR3996039}
for further geometric properties of this object).\medskip

Having settled questions~(Q1) and~(Q2) permits us to give an answer to question~(Q3)
as well, by exploiting a sliding method.
Indeed, as proved in~\cite{2019-1}, we have that the
stickiness phenomenon is
``generic", in the sense that
any small perturbation of any exterior datum is sufficient to produce boundary
discontinuities (hence, boundary continuity of planar nonlocal minimal graphs
should be considered as an ``exception'' to the ``typical''
case in which boundary jumps occur).
The precise statement of this result
is the following:

\begin{theorem}[Genericity of the stickiness phenomenon]\label{theorem-stickiness-is-generic}
	Let $u$ be an $s$-minimal graph in $(0,1)\times \R$ with smooth external
datum $u_0$. Suppose that
	\[u_0(0)=0=\lim_{x_1\searrow0}u(x_1).\]
Let~$\varphi\in C^{\infty}_0((-2,1),[0,+\infty))$ be not identically zero.
For every~$t>0$, let~$u^{(t)}$ be the $s$-minimal graph in $(0,1)\times \R$ with external
datum $u_0+t\varphi$. 
Then,
\begin{equation}\label{SR24}
\lim_{x_1\searrow 0}u^{(t)}(x_1)>0.\end{equation}
\end{theorem}

We observe that~\eqref{SR24} says that~$u^{(t)}$ always presents
the stickiness phenomenon for all~$t>0$, being the case~$t=0$
the only possible exception, namely a small positive bump always pushes up
the planar nonlocal minimal graphs in a discontinuous way at the boundary.
In other words,
if the ``unperturbed'' minimizer
is not sticky, then any positive, small and smooth perturbation of the datum will yield stickiness.
In this sense, our answer to question~(Q3)
is that stickiness is indeed quite a ``generic''
phenomenon representing the ``typical'' boundary behavior of nonlocal
minimal surfaces (with no counterpart in the theory of classical minimal surfaces).\medskip

    \begin{figure}[H]
    	\centering
    	\includegraphics[width=7.3cm]{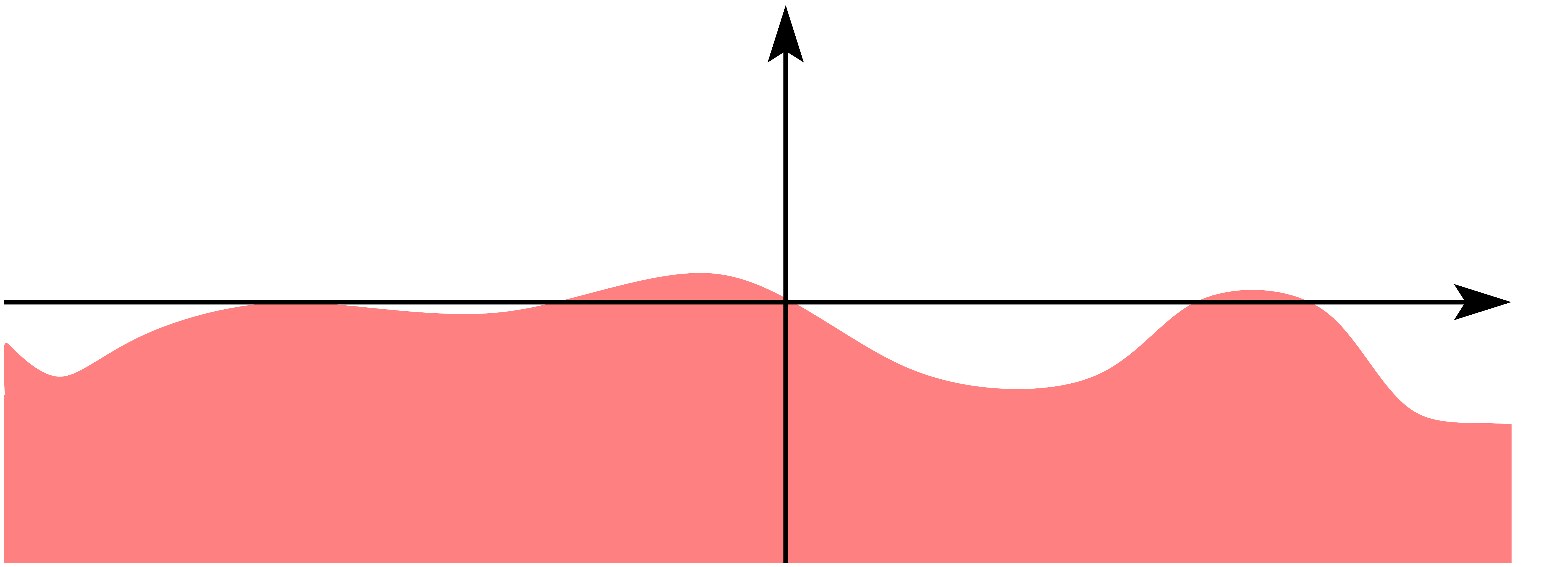} $\qquad$
    	\includegraphics[width=7.3cm]{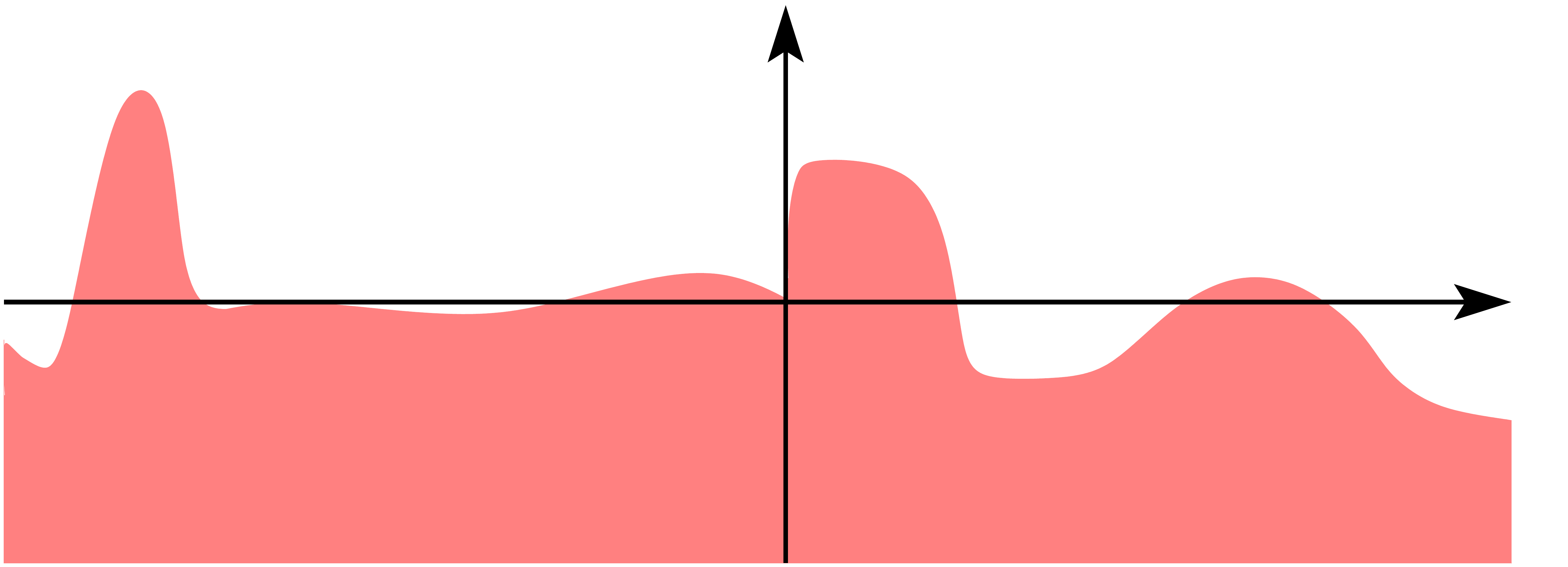}  
    	\caption{\textit {{The ``butterfly effect''
for the derivative of planar nonlocal minimal graphs (the additional
external bump on the left in the second picture can be taken arbitrarily
small, and arbitrarily far, still making the derivative infinite at the origin, here we have ``magnified''
this bump to improve the visibility effect).
}}}
    	\label{FIGBUT}
    \end{figure}

It is also interesting to observe that~\eqref{sh744}
and Theorem~\ref{FUN}, combined to Theorem \ref{theorem-stickiness-is-generic},
showcase a remarkable ``butterfly
effect'' for the derivative of planar nonlocal minimal graphs:
namely, if~$u_0'(0^-)=\ell$, for some~$\ell\in\R$,
and there is no stickiness effect at the origin, then also~$u'(0^+)=\ell$;
but as soon as a small, and possibly faraway, bump is placed
somewhere in the exterior datum, then suddenly~$|u'(0^+)|=+\infty$,
see Figure~\ref{FIGBUT}.
In this sense, the stickiness phenomenon also produces
generically the sudden divergence of the boundary derivative.
\medskip

The proof of Theorem \ref{theorem-stickiness-is-generic}
relies on a vertical sliding method. Specifically, one slides~$u$ and then
moves it till $u$ touches $u^{(t)}$ at some point.
Then, one reaches a contradiction using the Euler-Lagrange equation associated
with the $s$-perimeter minimization: for this,
it is however crucial to know that the Euler-Lagrange equation is
indeed satisfied at any point, and this is exactly the step in which
Theorem~\ref{VALID} comes into play.\medskip

Summarizing, question~(Q1) concerning
the boundary regularity of planar nonlocal minimal graphs is addressed
in~\cite{MR3532394} for discontinuous graphs and in Theorem~\ref{FUN}
for continuous graphs, thus leading to a general statement, valid both
for continuous and discontinuous graphs, as given in Theorem~\ref{QQ1},
saying that the boundary of planar nonlocal minimal graphs
is always a $C^{1,\frac{1+s}2}$-curve up to the boundary of the domain
(in a geometric sense). This in turn provides an answer for question~(Q2),
as in Theorem~\ref{VALID}, which ensures the validity of the Euler-Lagrange
equation at any point of the domain (including points at the boundary
of the domain, both in the case of continuous and discontinuous graphs).
This fact then allows one to exploit sliding methods, proving the genericity
of the stickiness phenomenon, thus answering question~(Q3)
as in Theorem~\ref{theorem-stickiness-is-generic}.

To complete the picture,
we now provide a sketch of the proof of Theorem \ref{FUN},
from which all the other results heavily depend.

\section*{Sketch of the proof of Theorem \ref{FUN}}

For simplicity, let us suppose that $u_0$ is zero on a left neighborhood of
the origin, say 
	\begin{equation}\label{simplicity-assumption}
	u_0(x_1)= 0\quad \mbox{ for every }x_1\in[-h,0],
	\end{equation} 
for some~$h>0$. We stress that~\eqref{simplicity-assumption}
is a slightly simplifying assumption
when compared to assumption~\eqref{u0reg}, but the arguments
presented here
would carry over, up to technical complications, just assuming that~$u_0$
is sufficiently smooth in a left neighborhood of the origin (full details available
in~\cite{2019-1}).

	Now, roughly speaking, the idea of the proof is to
``look at the worst possible scenarios'' and
``rule out all the other possibilities''.

\begin{figure}[H]
    	\centering
    	\includegraphics[width=6cm]{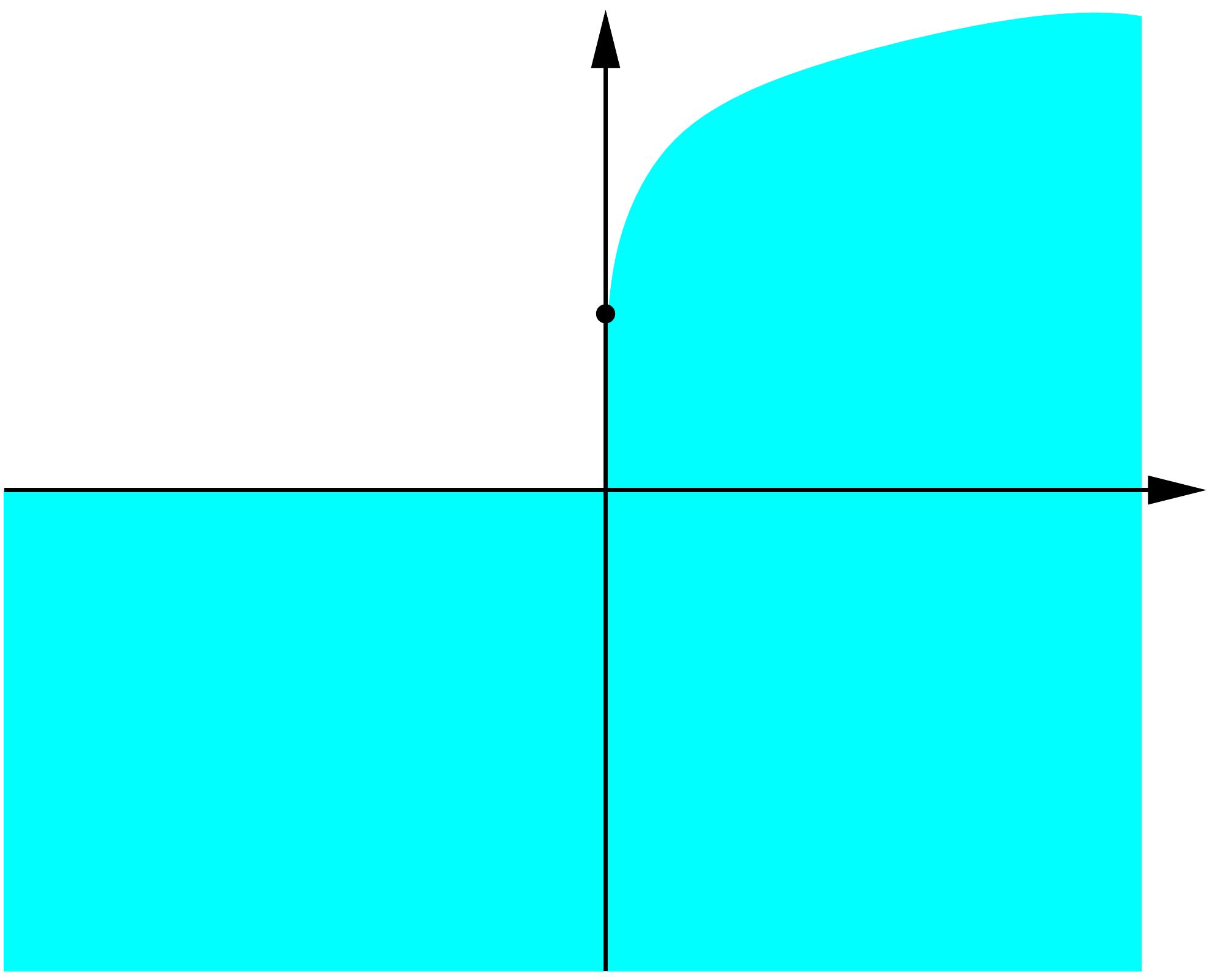} $\qquad$
    	\includegraphics[width=6cm]{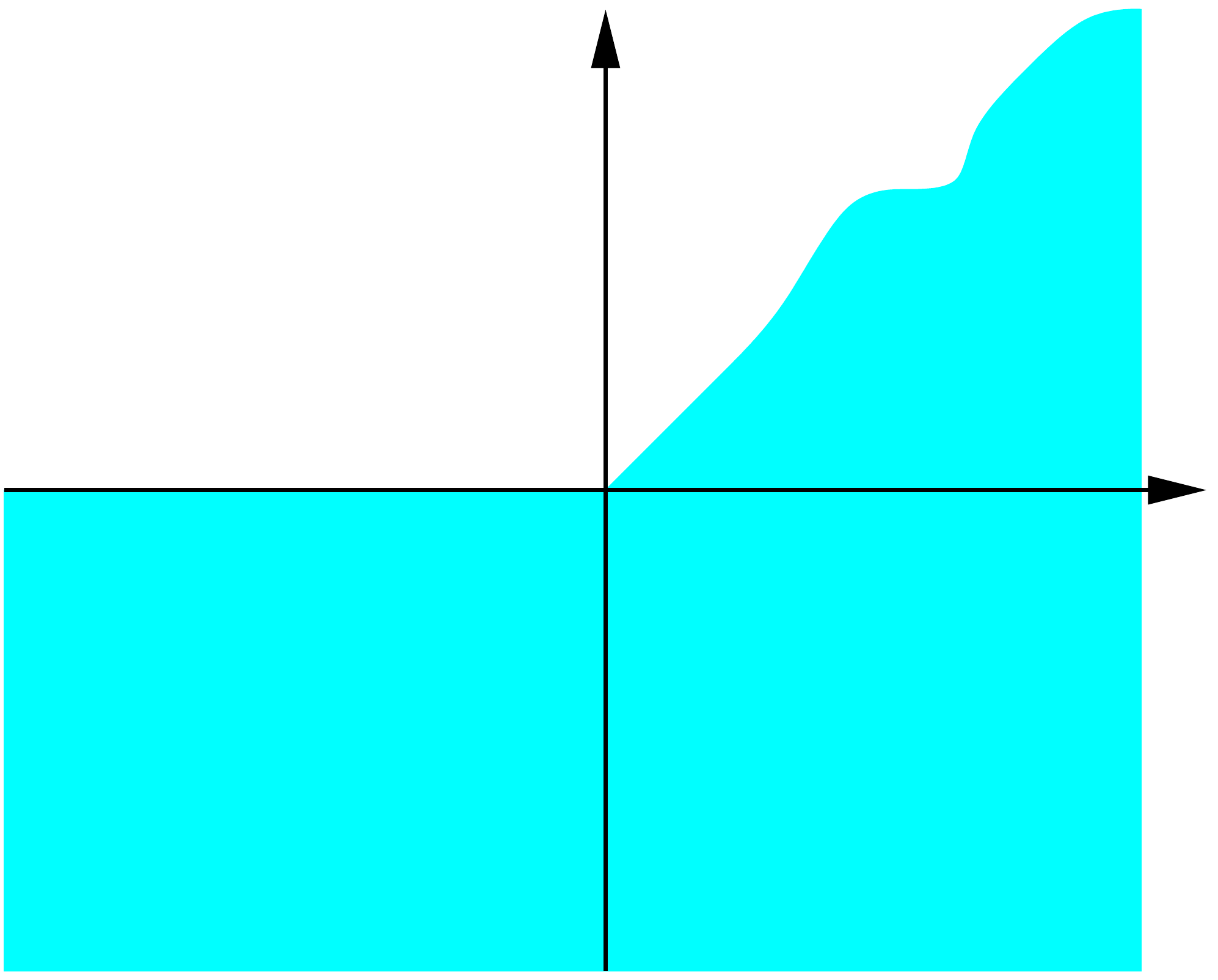} \\ $\qquad$ \\
    	\includegraphics[width=6cm]{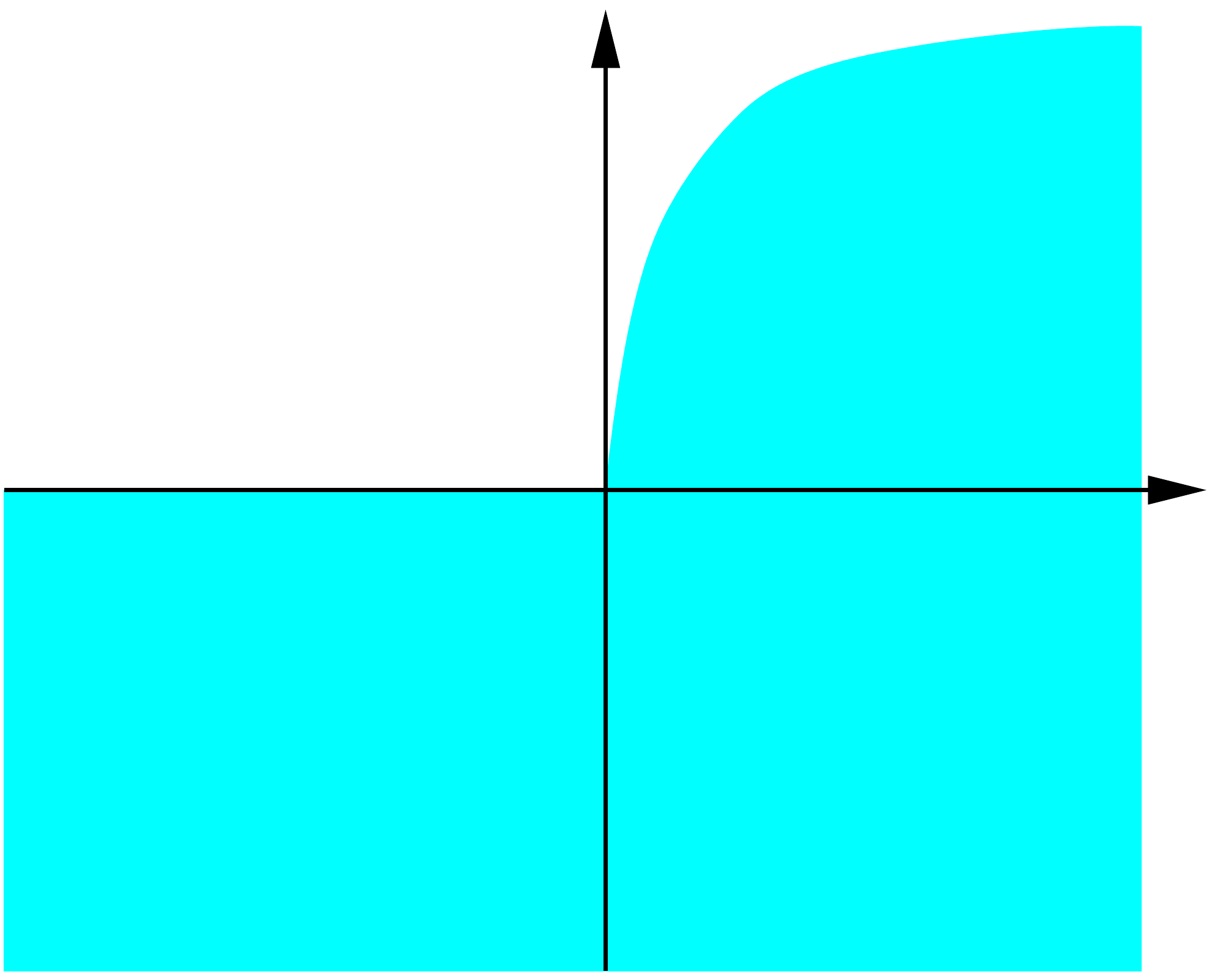} $\qquad$
    	\includegraphics[width=6cm]{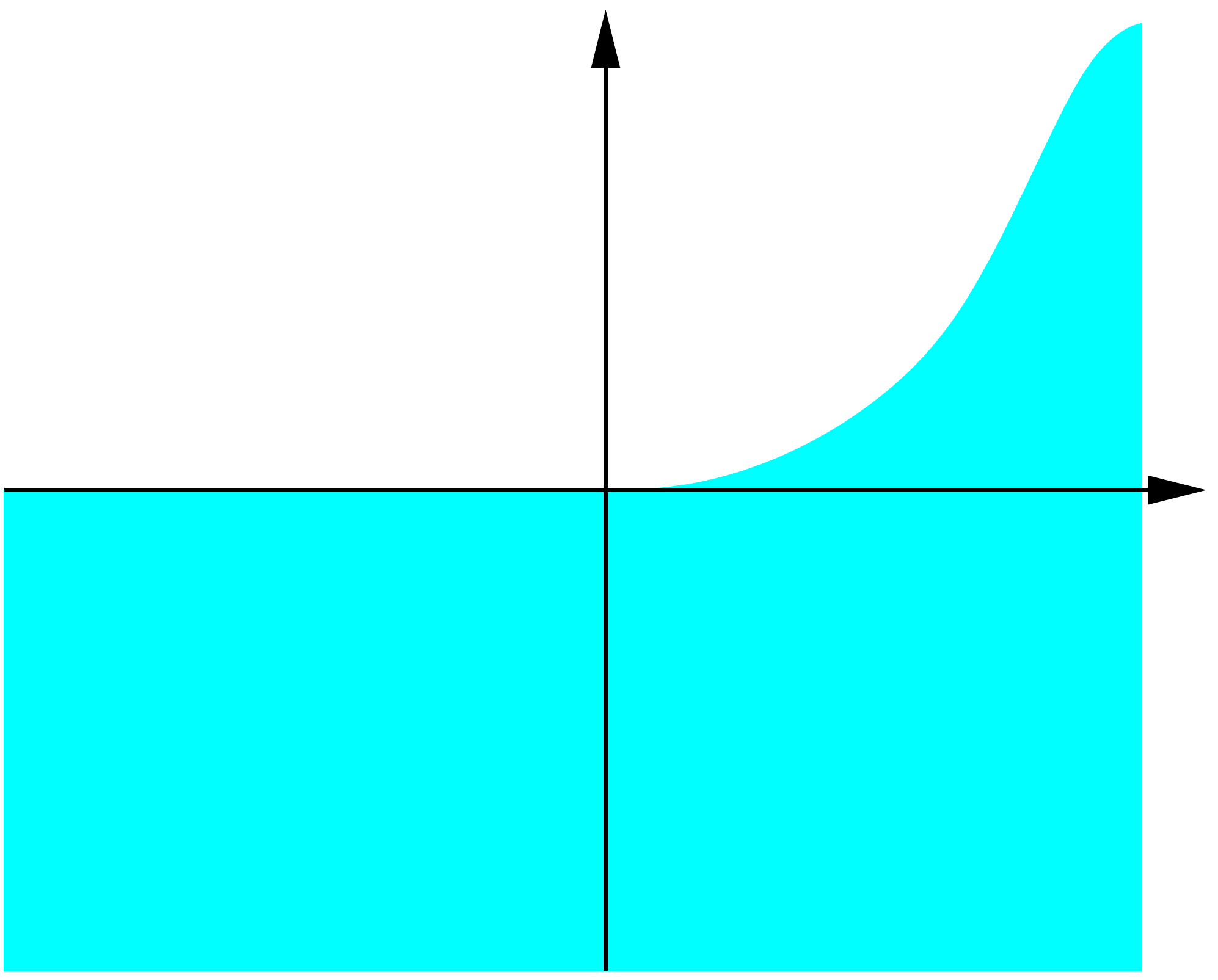}
    	\caption{\textit {{What are the possible boundary behaviors of nonlocal
    				minimal graphs in $(0,1)$
    				vanishing in $(-1,0]$?}}}
    	\label{FIG6}
    \end{figure}

To make this strategy concrete,
we can consider the 
prototype situations
embodied by the
following\footnote{We remark that these cases
do not really exhaust all the possibilities, but they nevertheless provide
a very good indication of what's going on in the general situation.
For full details on the proof of Theorem \ref{FUN},
we refer to~\cite{2019-1}.} cases (see
Figure~\ref{FIG6}):
	\begin{enumerate}
		\item[(i)] $u$ has a jump discontinuity at the origin,
thus exhibiting the stickiness phenomenon -- but
this occurrence is ruled out in this case by assumption~\eqref{continuity-condition};
		\item[(ii)] \label{Lipschitz-possibility}$u$ is Lipschitz continuous in $[-h,1/2],$ but not better than this;
		\item[(iii)] \label{Holder-possibility}$u\in C^{\alpha}([-h,1/2]),$ for some~$0<\alpha<1,$ but not better than this;
		\item[(iv)] \label{C-1-possibility}$u\in C^1([-h,1/2]),$ but~$u\not\in C^{1,\gamma}([-h,1/2]).$
	\end{enumerate}
Hence, to convince ourselves of the validity of
Theorem~\ref{FUN}, it is necessary to exclude
the possibilities described in~(ii), (iii) and~(iv)
(and also to obtain a uniform bound on the H\"older exponent 
of the derivative of~$u$).

To do so, it is convenient to consider the blow-up
limits corresponding to~(ii), (iii) and~(iv)
and try to understand their relations
with the original picture.

If $E_u$ is as in \eqref{definition-E-u}, for the sake of shortness
we denote it by $E$, and we define the blow-up sequence $E_k$ of $E,$
with $k\in \mathbb{N},$ as
    \begin{equation}\label{definition-blow-up-sequence}
    E_k:=kE=\left\{k(x_1,x_2),\hspace{0.2cm}(x_1,x_2)\in E\right\}.
    \end{equation} 
As a technical remark, we recall that
the existence of the blow-up limit, 
that is the limit
as~$k\to+\infty$, possibly up to a subsequence, of the set in~\eqref{definition-blow-up-sequence},
typically follows from suitable density estimates
(in this framework, these density
estimates need to be centered at a boundary point, 
and the setting in~\eqref{simplicity-assumption} allows one to
extend the interior estimates to the case under consideration,
see Lemma 2.1 in \cite{2019-1} for details).
\medskip

Now, we would like to say that the blow-up limit is a cone.
This usually relies on a specific monotonicity formula, and, in our framework,
such a precise monotonicity formula is not available.
To circumvent this difficulty, it is convenient to replace
the previous blow-up limit with a second blow-up limit
(that is, one considers a blow-up sequence obtained from
the first blow-up limit, and then takes the limit of this
new sequence).
The advantage of this second blow-up procedure is that 
the first blow-up limit is already a halfplane in~$\{x_1<0\}$,
thanks to~\eqref{simplicity-assumption}; consequently,
every element of the new blow-up sequence is already a halfplane in~$\{x_1<0\}$.
{F}rom this, the proof of the monotonicity formula in~\cite{MR2675483}
carries over for the second blow-up sequence, thus ensuring that
the second blow-up limit is indeed a cone (full details
of this construction can be found in Lemma 2.2 in \cite{2019-1}).\medskip

\begin{figure}[H]
	\centering
	\includegraphics[width=6cm]{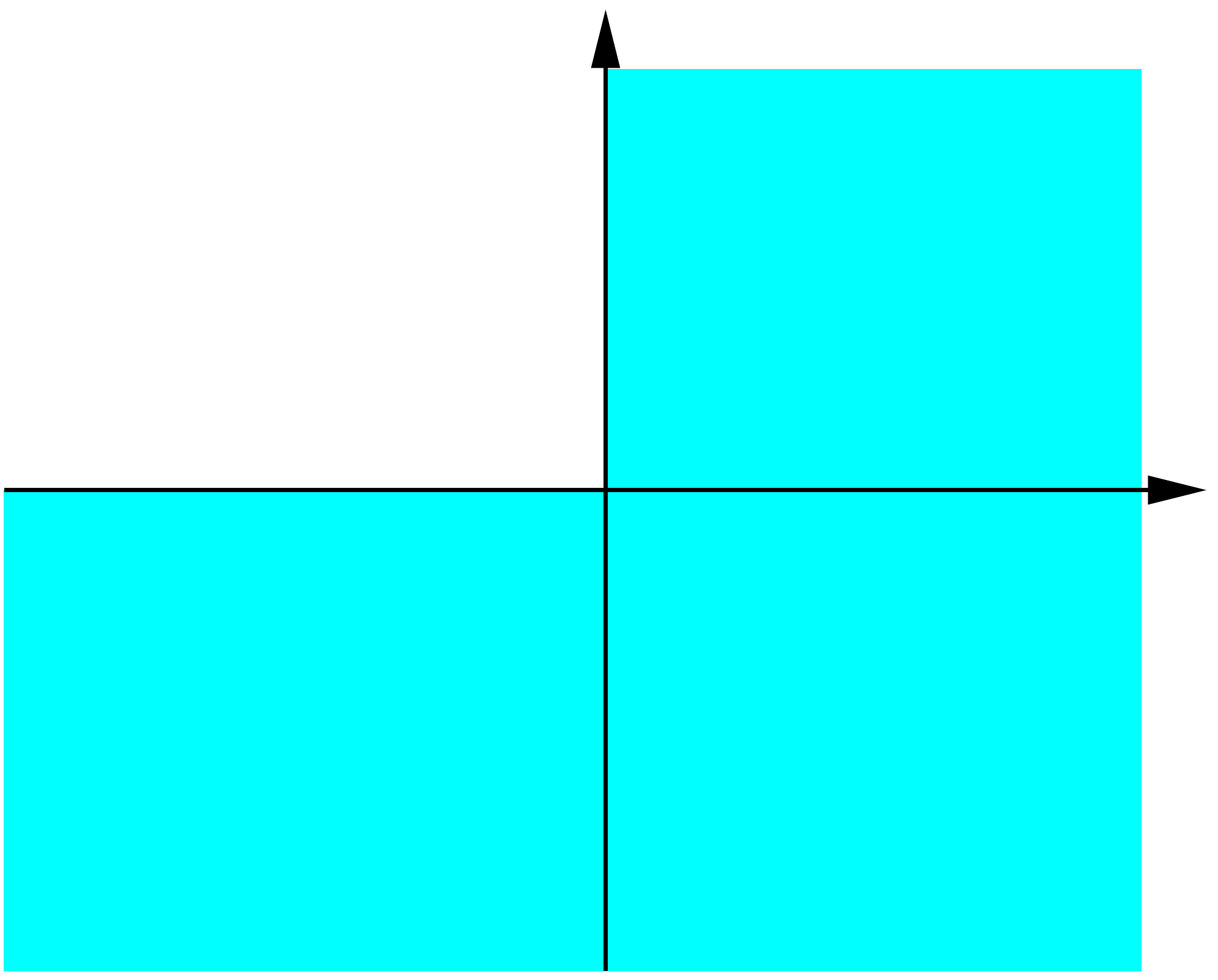} $\qquad$
	\includegraphics[width=6cm]{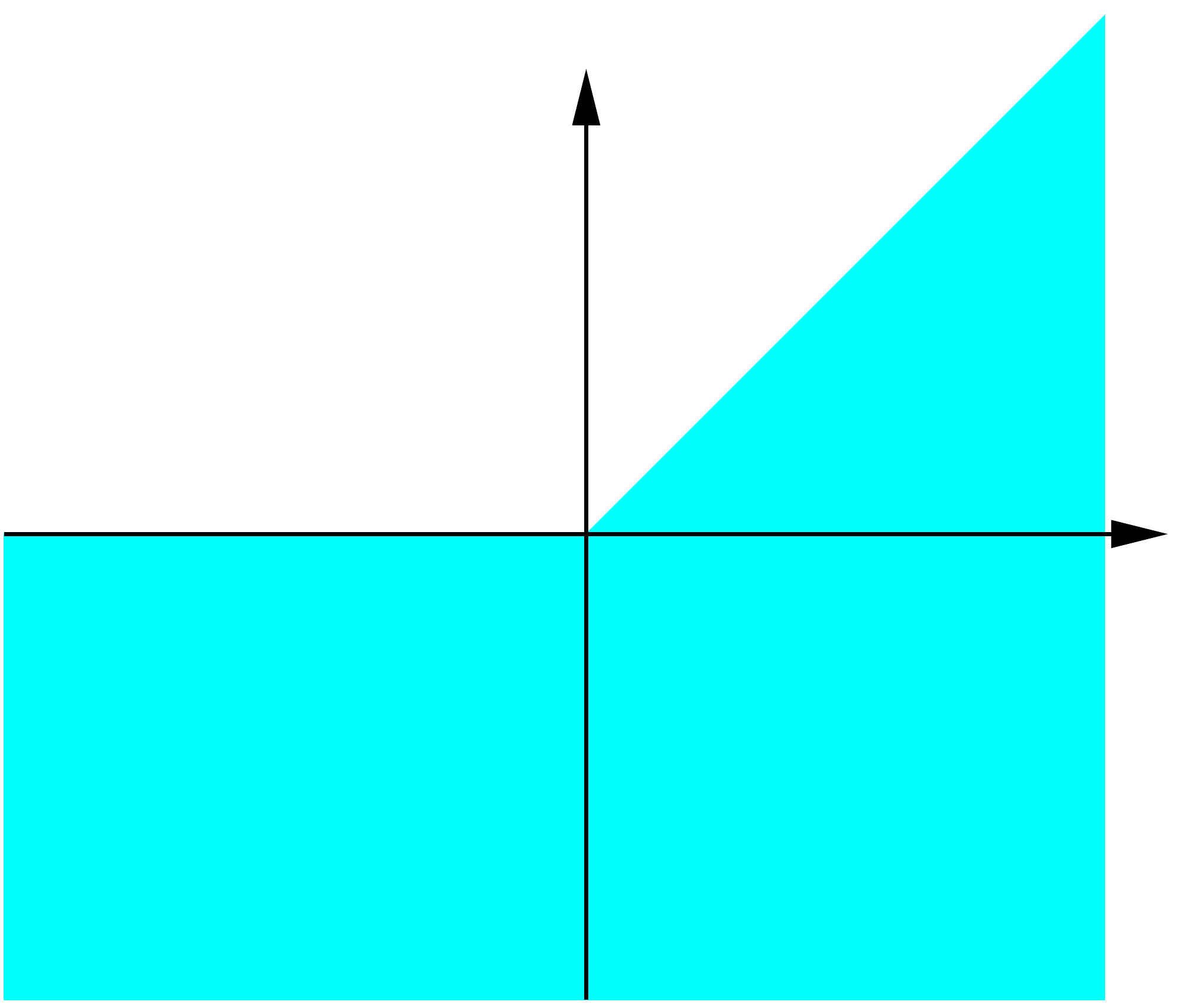} \\ $\qquad$ \\
	\includegraphics[width=6cm]{101c.pdf} $\qquad$
	\includegraphics[width=6cm]{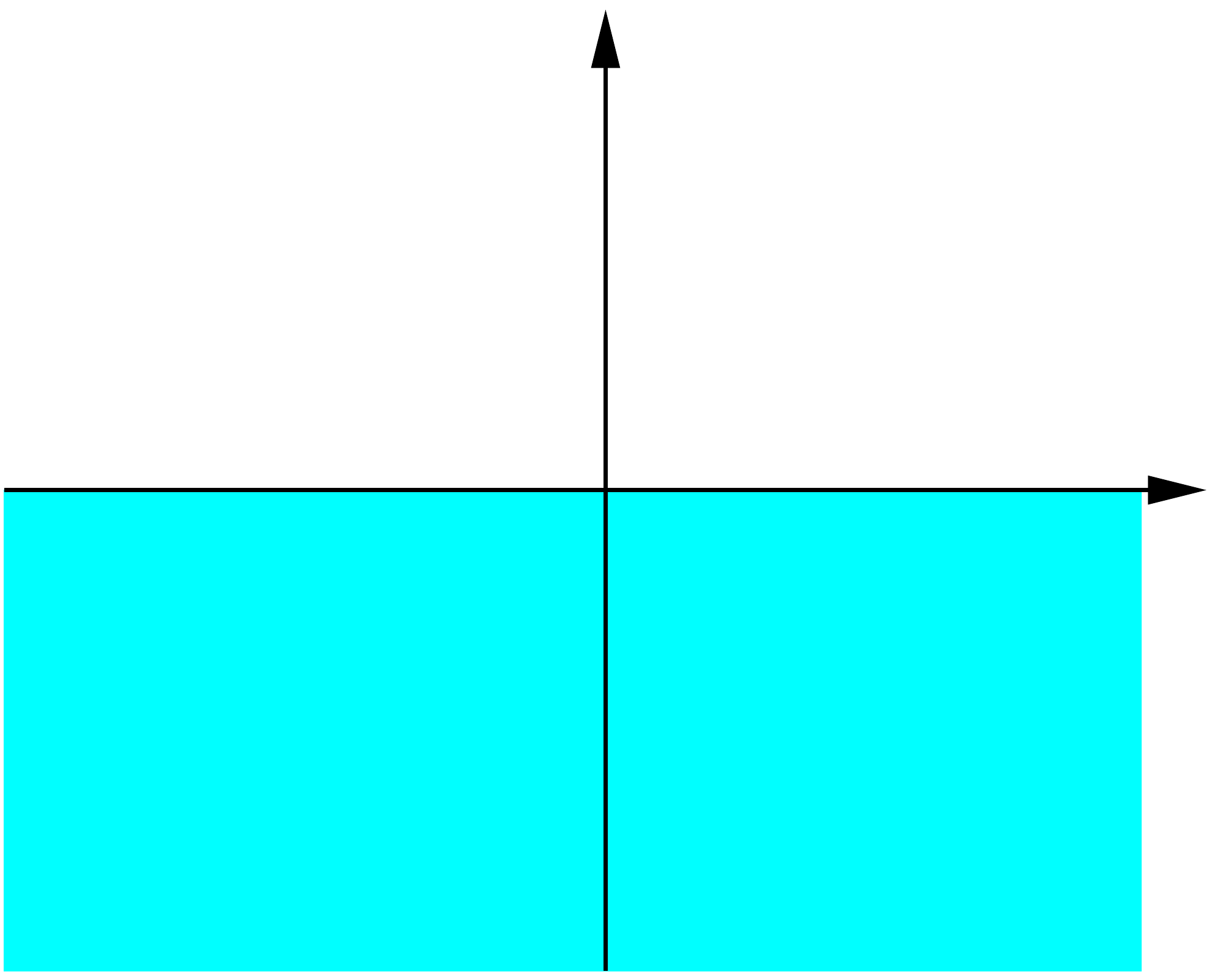}
	\caption{\textit {{Blow-up's of the four possibilities depicted in Figure \ref{FIG6}.}}}
	\label{FIG7}
\end{figure}

We denote the second blow-up limit by~$E_{00}$, and we recall that,
in view of~\eqref{simplicity-assumption}, we have that 
    \[
    E_{00}\cap \left\{x_1<0\right\}=\left\{x_1<0,\quad x_2<0\right\}.
    \]
See also Figure \ref{FIG7} for a description of
the second blow-up limits corresponding to the possibilities
depicted in Figure~\ref{FIG6}.
Comparing with the possibilities~(ii), (iii) and~(iv), that should be ruled
out in order to establish Theorem~\ref{FUN},
we obtain the following scenarios for the second blow-up's:
    \begin{enumerate}
    	\item[(ii)'] \label{Lipschitz-blow-up} $E_{00}\cap\left\{x_1
>0\right\}=\left\{x_2<bx_1\right\}\cap\left\{x_1>0\right\}$,
for some~$b\in\R$, which is the second blow-up limit corresponding to possibility~(ii);
    	\item[(iii)'] \label{Holder-blow-up}$E_{00}\cap \left\{x_1>0\right\}=\left\{x_1>0\right\},$ which is the second blow-up limit
corresponding to
alternative~(iii);
    	\item[(iv)'] \label{C-1-blow-up}$E_{00}=\left\{x_2<0\right\},$ that is $E_{00}$ is a half-plane, which is the second blow-up limit corresponding to
possibility~(iv).
    \end{enumerate}
Hence, our sketch of the proof of Theorem~\ref{FUN} would be completed
once we eliminate the possibilities in~(ii)', (iii)' and~(iv)'.
\medskip

    First, we proceed to exclude possibility~(ii)'. For this, using the minimality of~$
E_{00}$ in $\left\{x_1>0\right\}$
we have that
\begin{equation}\label{U7-coneuy-1}\int\limits_{\R^2}\frac{\chi_{\R^2\setminus
E_{00}}(y)-\chi_{ E_{00}}(y)}{
\left|p-y\right|^{2+s}}\,dy=0,\end{equation}
where~$p:=(1,b)$.

On the other hand, we see that
\begin{equation}\label{U7-coneuy-2}
\int\limits_{\R^2}\frac{\chi_{\R^2\setminus E_{00}}(y)-\chi_{E_{00}}(y
)}{\left|p-y\right|^{2+s}}\,dy\ne 0,\end{equation}
since the contribution of the set and the one of its complement
do not cancel each other
(compare the symmetric regions arising after drawing the tangent line
passing through~$p$). 
\medskip

The contradiction arising from~\eqref{U7-coneuy-1}
and~\eqref{U7-coneuy-2} rules out possibility~(ii)', and we now want
to exclude possibility~(iii)'. We observe that, to exclude
this possibility, one cannot only rely on blow-up type analysis,
since the same blow-up limit as the one in~(iii)'
is also achieved when~$u$ 
is discontinuous: that is, possibility~(i)
would produce the same blow-up picture as possibility~(iii),
but the original nonlocal minimal graphs present obvious structural differences.
Therefore, the strategy to eliminate~(iii)'
has to take into account the original sets with a finer analysis,
and indeed we aim at showing that possibility~(iii)' can only come
from discontinuous nonlocal minimal graph~$u$ (and this possibility,
corresponding to~(i), was already ruled out in light of~\eqref{continuity-condition}).\medskip

In this sense, the strategy to eliminate possibility~(iii)'
consists in proving that ``thick $s$-minimal sets are necessarily full''
(or, equivalently, considering complement sets,
``narrow $s$-minimal sets are necessarily void"). The precise statement,
which is a particular case of Proposition~3.1 in~\cite{2019-1},
goes as follows:

    \begin{proposition}\label{proposition-sliding-method}
    	Let $\lambda>0.$ There exist $M_0>1$ and $\mu_0\in (0,1)$ such that
if $M\geq M_0$ and $\mu \in (0,\mu_0]$ the following claim holds true.

    	Let $F\subset \R^2$ be $s$-minimal in $
(0,M)\times (-4,4)$.
If 
    	\begin{align}\label{24637rtegub293yreghf}
    	&F\cap \left\{x_1 \in (-M,0)\right\}\;= \;\left\{x_2\le0\right\},\\
    	\label{E-in-a-small-slab-hypothesis}\mbox{and}\qquad& \Big(
(0,M)\times (-M,M)\Big)\setminus F\;\subseteq \;\left\{x_1\in (0,\mu)\right\},
    	\end{align}
    	then
    	\begin{equation}\label{rq6ewt9238eu} \left(0,\frac{M}2\right)\times (-1,1)\;\subseteq \;F. \end{equation}
    \end{proposition}

The idea to prove Proposition~\ref{proposition-sliding-method}
is to argue by contradiction exploiting a sliding method.
Namely, 
if the thesis in~\eqref{rq6ewt9238eu} were false,
one could take a ball inside $F$ and slide it till it touches the complement of~$
F$ at some point~$q$. In this framework, one obtains the existence of
a ball~$B\subseteq F$, with~$q\in \partial B\cap \partial F$.
The strategy is to show that the $s$-mean curvature of $F$ at $q$
is strictly negative, thus
contradicting the minimality of $F$.\medskip

\begin{figure}[H]
	\centering
	\includegraphics[height=8cm]{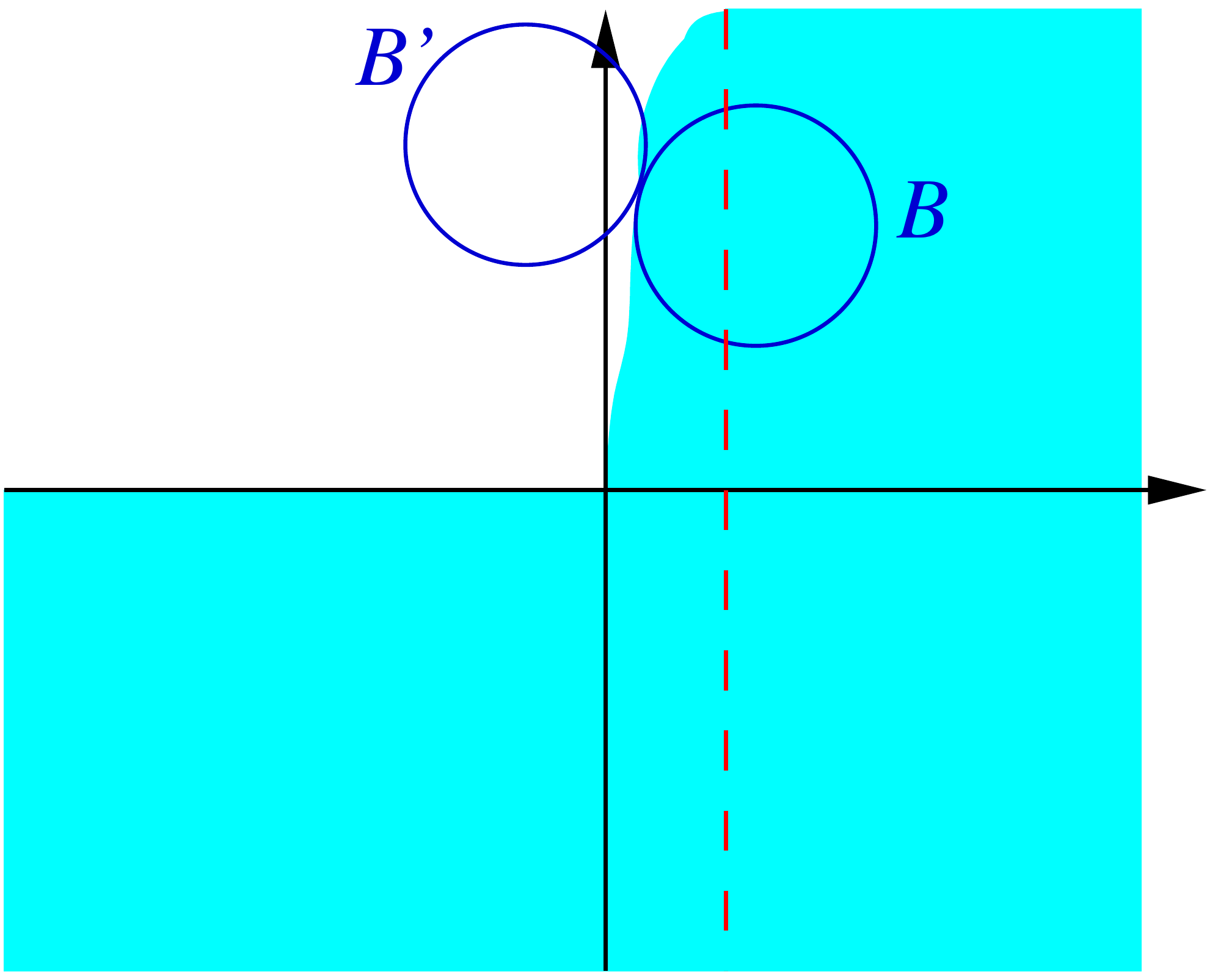}
	\caption{\textit{{Using symmetric balls to detect integral
cancellations.}}}
	\label{XCdFIG9}
\end{figure}

To compute the $s$-mean curvature of $F$ in $q$,
it is convenient to consider the symmetric ball to $B$ with respect to
the tangent plane 
through $q$ and denote it by $B'$,
see Figure~\ref{XCdFIG9}.
By hypothesis \eqref{E-in-a-small-slab-hypothesis}, we know that the complement
of~$F$ (in the domain~$\{x_1>0\}$, up to a large cylinder)
is contained in a small slab near the vertical axis.
Hence, the integral contributions for the $s$-mean
curvature (recall the left hand side in~\eqref{EELLQQ})
are ``mostly negative'', coming predominantly from
points in the set~$F$, with the possible exception of
the points in the complement lying in the small slab~$\left\{x_1\in (0,\mu)\right\}$.
Near the contact point~$q$, the positive contributions
coming from these points are ``negligible'' as long as~$\mu$
is sufficiently small, since the singularity of the kernel
is compensated by the symmetric integration over the balls~$B$ and~$B'$,
with the full ball~$B$ providing negative contributions.
Similarly, far from~$q$, the negative terms coming from the set~$F$
provide a negligible contribution
to the $s$-mean curvature, since the singularity of the kernel
plays little role away from~$q$, and the weighted measure
of the narrow slab is small with~$\mu$.\medskip

These quantitative arguments
establish
Proposition~\ref{proposition-sliding-method} 
(see again Proposition~3.1 in~\cite{2019-1} for full details).
With this, one can rule out possibility~(iii)'
by arguing as follows. {F}rom~(iii)', one knows that the
blow-up limit in~$\{x_1>0\}$ is either full
or void. Let us consider the first case (up to changing a set with its
complement), namely suppose that
\begin{equation}\label{gsdfP} E_{00}\cap\{x_1>0\}=
(0,+\infty)\times\R.\end{equation}
Then, up to a subsequence, a suitable blow-up sequence~$E_k$
(recall~\eqref{definition-blow-up-sequence}) would lie locally in the vicinity
of~$E_{00}$ for a suitably large~$k$. {F}rom this and~\eqref{gsdfP},
one sees that~$E_k$ satisfies~\eqref{E-in-a-small-slab-hypothesis}
for~$k$ sufficiently large (possibly depending on the thresholds~$M_0$
and~$\mu_0$ in Proposition~\ref{proposition-sliding-method},
and notice also that in this setting~$E_k$
satisfies~\eqref{24637rtegub293yreghf}
as a consequence of~\eqref{simplicity-assumption}).
Then, one can apply
Proposition~\ref{proposition-sliding-method} to~$F:=E_k$. Consequently,
from~\eqref{24637rtegub293yreghf}
and~\eqref{rq6ewt9238eu}, one obtains that the graph describing~$E_k$
has a jump discontinuity at the origin.
Scaling back, this gives that~$u$ has a jump discontinuity at the origin.
This is in contradiction with~\eqref{continuity-condition},
and therefore possibility~(iii)' (and hence~(iii)) has been excluded.
\medskip

It remains to rule out possibility~(iv)' (and hence~(iv)).
To this end, we need to prove that once the blow-up limit is a half-plane,
then necessarily the original $s$-minimal graph
was already differentiable at the origin, with a precise estimate on the
H\"older exponent of the derivative (we stress that
controlling the H\"older exponent of the derivative is a crucial step
in order to deduce the results in Theorems~\ref{QQ1}, \ref{VALID}
and~\ref{theorem-stickiness-is-generic}
from Theorem~\ref{FUN}).

The idea of the proof now
consists in using ``vertical rescalings'' for an ``improvement of flatness''
(once we know that the solution is sufficiently flat
at a large scale,
then it is necessarily even flatter at a smaller scale).
Differently than other improvement of flatness methods,
which were designed in the interior of the domain (see~\cite{MR2675483}),
our setting requires us to achieve this enhanced regularity at boundary points.
To this aim, one considers vertical rescalings
and proves convergence to some function~$\bar{u}$,
which satisfies~$(-\Delta)^{\sigma}\bar{u}=0$ in~$(0,+\infty)$,
with~$\sigma=\frac{1+s}{2}$, and~$\bar u=0$ in~$(-\infty,0)$.
The linear theory of fractional equations (see e.g.~\cite{MR3168912})
would only ensure that~$\bar u$ is H\"older continuous at the origin,
but our objective is to prove that in fact~$\bar u$ is more regular, thus
producing the desired enhanced regularity for the original function~$u$,
by bootstrapping such improvement of flatness method.\medskip

\begin{figure}[H]
	\centering
	\includegraphics[height=8cm,width=15cm]{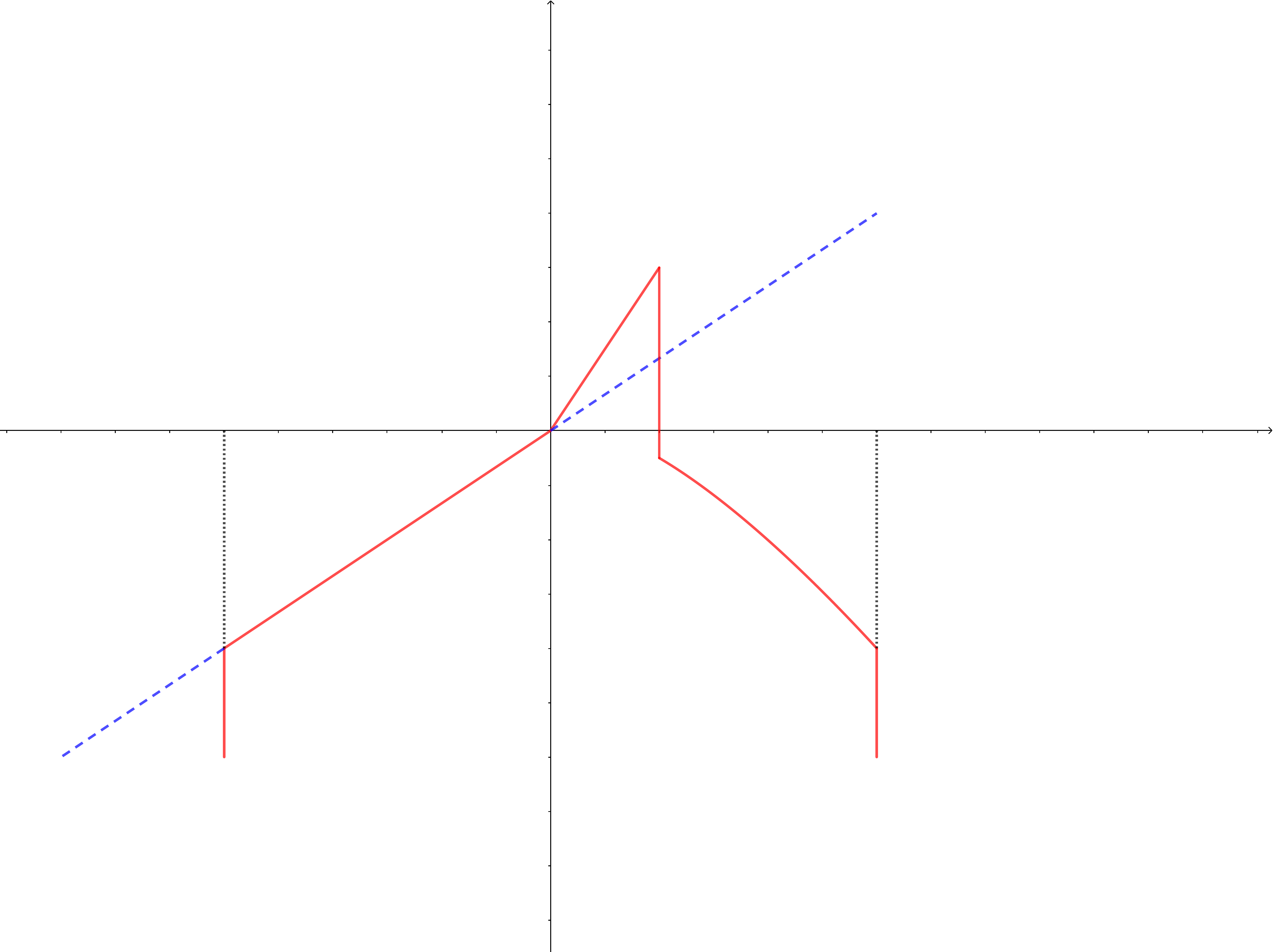}
	\caption{\textit{{Shape of the corner-like barrier.}}}
	\label{FIG9}
\end{figure}

Concretely, one deduces from the linear theory of fractional equations
that, for small~$x_1>0$,
\begin{equation}\label{hdbc12345}\bar{u}(x_1)=a_0 x_1^{\sigma}+O(x_1^{\sigma+1})\end{equation}
for a suitable~$a_0\in \R$. Our goal is to show that
\begin{equation}\label{hdbc123452}
a_0=0,\end{equation}
thus improving the boundary regularity in this specific case. To do this, we construct a
suitable corner-like
barrier (see Figure~\ref{FIG9} here, and Lemma 7.1 in \cite{2019-1} for full details).

Roughly speaking, one can juggle parameters to make
the subgraph depicted in Figure~\ref{FIG9} have negative fractional
mean curvature
in the vicinity of the origin.
Intuitively, this is possible thanks to a ``purely nonlocal effect'': indeed,
in the classical case, the segment near the origin in Figure~\ref{FIG9}
would produce a zero curvature (thus making the argument invalid
for classical minimal surfaces), while in the nonlocal setting the concave
corner at the origin produces a very negative fractional curvature (actually,
equal to~$-\infty$ at the origin). This negative contribution survives
after the bending of the barrier at the side of Figure~\ref{FIG9}
(which is needed in order to place the barrier ``below the solution
at infinity'').

Then, to prove~\eqref{hdbc123452}, one argues by contradiction,
supposing, for instance, that~$a_0>0$. Then, using~\eqref{hdbc12345},
one sees that the barrier in Figure~\ref{FIG9} can be slided from
below the original $s$-minimal graph~$u$. By maximum principle
(and noticing the linear growth of the barrier
in Figure~\ref{FIG9} for~$x_1>0$ small),
this gives that~$u$ lies above a linear function
for~$x_1>0$ small. Consequently, the corresponding blow-up limit would be
as in possibility~(ii)'. Since this alternative has been already ruled out, we obtain
a contradiction, thus establishing~\eqref{hdbc123452}.\medskip

{F}rom~\eqref{hdbc123452}, the improvement of flatness method
kicks in, thus producing the desired enhancement regularity result
that rules out the last possibility, finally leading to the completion of
the proof of Theorem~\ref{FUN}.

\vfill

{\footnotesize

\noindent {\em Addresses:} \\

Serena Dipierro \& Enrico Valdinoci.
Department of Mathematics
and Statistics,
University of Western Australia,
35 Stirling Hwy, Crawley WA 6009, Australia.\\

Aleksandr Dzhugan \& Nicol\`o Forcillo.
Dipartimento di Matematica,
Alma Mater Studiorum Universit\`a di Bologna,
Piazza di Porta San Donato 5, 40126 Bologna, Italy.\\

\smallskip

\noindent{\em Emails:}\\

serena.dipierro@uwa.edu.au \\

aleksandr.dzhugan2@unibo.it \\

nicolo.forcillo2@unibo.it \\

enrico.valdinoci@uwa.edu.au \\  
}
\begin{bibdiv}
\begin{biblist}

\bib{MR3230079}{article}{
   author={Abatangelo, Nicola},
   author={Valdinoci, Enrico},
   title={A notion of nonlocal curvature},
   journal={Numer. Funct. Anal. Optim.},
   volume={35},
   date={2014},
   number={7-9},
   pages={793--815},
   issn={0163-0563},
   review={\MR{3230079}},
   doi={10.1080/01630563.2014.901837},
}

\bib{MR2765717}{article}{
   author={Ambrosio, Luigi},
   author={De Philippis, Guido},
   author={Martinazzi, Luca},
   title={Gamma-convergence of nonlocal perimeter functionals},
   journal={Manuscripta Math.},
   volume={134},
   date={2011},
   number={3-4},
   pages={377--403},
   issn={0025-2611},
   review={\MR{2765717}},
   doi={10.1007/s00229-010-0399-4},
}

\bib{MR3331523}{article}{
   author={Barrios, Bego\~{n}a},
   author={Figalli, Alessio},
   author={Valdinoci, Enrico},
   title={Bootstrap regularity for integro-differential operators and its
   application to nonlocal minimal surfaces},
   journal={Ann. Sc. Norm. Super. Pisa Cl. Sci. (5)},
   volume={13},
   date={2014},
   number={3},
   pages={609--639},
   issn={0391-173X},
   review={\MR{3331523}},
}

\bib{MR1942116}{article}{
   author={Brezis, Kh.},
   title={How to recognize constant functions. A connection with Sobolev
   spaces},
   language={Russian, with Russian summary},
   journal={Uspekhi Mat. Nauk},
   volume={57},
   date={2002},
   number={4(346)},
   pages={59--74},
   issn={0042-1316},
   translation={
      journal={Russian Math. Surveys},
      volume={57},
      date={2002},
      number={4},
      pages={693--708},
      issn={0036-0279},
   },
   review={\MR{1942116}},
   doi={10.1070/RM2002v057n04ABEH000533},
}

\bib{MR1945278}{article}{
   author={Bourgain, Jean},
   author={Brezis, Ha\"{\i}m},
   author={Mironescu, Petru},
   title={Limiting embedding theorems for $W^{s,p}$ when $s\uparrow1$ and
   applications},
   note={Dedicated to the memory of Thomas H. Wolff},
   journal={J. Anal. Math.},
   volume={87},
   date={2002},
   pages={77--101},
   issn={0021-7670},
   review={\MR{1945278}},
   doi={10.1007/BF02868470},
}

\bib{CLALUC}{article}{
author={Bucur, Claudia},
author={Lombardini, Luca},
title={Asymptotics as $s \searrow 0$ of the nonlocal nonparametric Plateau problem with obstacles},
      journal = {In preparation},
}

\bib{MR3926519}{article}{
   author={Bucur, Claudia},
   author={Lombardini, Luca},
   author={Valdinoci, Enrico},
   title={Complete stickiness of nonlocal minimal surfaces for small values
   of the fractional parameter},
   journal={Ann. Inst. H. Poincar\'{e} Anal. Non Lin\'{e}aire},
   volume={36},
   date={2019},
   number={3},
   pages={655--703},
   issn={0294-1449},
   review={\MR{3926519}},
   doi={10.1016/j.anihpc.2018.08.003},
}

\bib{CAB-C}{article}{
       author = {Cabr\'e, Xavier},
        title = {Calibrations and null-Lagrangians for nonlocal perimeters and an application to the viscosity theory},
      journal = {arXiv e-prints},
date={2019},
          eid = {arXiv:1905.10790},
        pages = {arXiv:1905.10790},
archivePrefix = {arXiv},
       eprint = {1905.10790},,
       adsurl = {https://ui.adsabs.harvard.edu/abs/2019arXiv190510790C},
}

\bib{MR3934589}{article}{
   author={Cabr\'{e}, Xavier},
   author={Cozzi, Matteo},
   title={A gradient estimate for nonlocal minimal graphs},
   journal={Duke Math. J.},
   volume={168},
   date={2019},
   number={5},
   pages={775--848},
   issn={0012-7094},
   review={\MR{3934589}},
   doi={10.1215/00127094-2018-0052},
}

\bib{MR3881478}{article}{
   author={Cabr\'{e}, Xavier},
   author={Fall, Mouhamed Moustapha},
   author={Sol\`a-Morales, Joan},
   author={Weth, Tobias},
   title={Curves and surfaces with constant nonlocal mean curvature: meeting
   Alexandrov and Delaunay},
   journal={J. Reine Angew. Math.},
   volume={745},
   date={2018},
   pages={253--280},
   issn={0075-4102},
   review={\MR{3881478}},
   doi={10.1515/crelle-2015-0117},
}

\bib{MR3744919}{article}{
   author={Cabr\'{e}, Xavier},
   author={Fall, Mouhamed Moustapha},
   author={Weth, Tobias},
   title={Delaunay hypersurfaces with constant nonlocal mean curvature},
   language={English, with English and French summaries},
   journal={J. Math. Pures Appl. (9)},
   volume={110},
   date={2018},
   pages={32--70},
   issn={0021-7824},
   review={\MR{3744919}},
   doi={10.1016/j.matpur.2017.07.005},
}

\bib{MR3770173}{article}{
   author={Cabr\'{e}, Xavier},
   author={Fall, Mouhamed Moustapha},
   author={Weth, Tobias},
   title={Near-sphere lattices with constant nonlocal mean curvature},
   journal={Math. Ann.},
   volume={370},
   date={2018},
   number={3-4},
   pages={1513--1569},
   issn={0025-5831},
   review={\MR{3770173}},
   doi={10.1007/s00208-017-1559-6},
}

\bib{MR3532394}{article}{
   author={Caffarelli, L.},
   author={De Silva, D.},
   author={Savin, O.},
   title={Obstacle-type problems for minimal surfaces},
   journal={Comm. Partial Differential Equations},
   volume={41},
   date={2016},
   number={8},
   pages={1303--1323},
   issn={0360-5302},
   review={\MR{3532394}},
   doi={10.1080/03605302.2016.1192646},
}

\bib{MR2675483}{article}{
   author={Caffarelli, L.},
   author={Roquejoffre, J.-M.},
   author={Savin, O.},
   title={Nonlocal minimal surfaces},
   journal={Comm. Pure Appl. Math.},
   volume={63},
   date={2010},
   number={9},
   pages={1111--1144},
   issn={0010-3640},
   review={\MR{2675483}},
   doi={10.1002/cpa.20331},
}

\bib{MR2564467}{article}{
   author={Caffarelli, Luis A.},
   author={Souganidis, Panagiotis E.},
   title={Convergence of nonlocal threshold dynamics approximations to front
   propagation},
   journal={Arch. Ration. Mech. Anal.},
   volume={195},
   date={2010},
   number={1},
   pages={1--23},
   issn={0003-9527},
   review={\MR{2564467}},
   doi={10.1007/s00205-008-0181-x},
}

\bib{MR2782803}{article}{
   author={Caffarelli, Luis},
   author={Valdinoci, Enrico},
   title={Uniform estimates and limiting arguments for nonlocal minimal
   surfaces},
   journal={Calc. Var. Partial Differential Equations},
   volume={41},
   date={2011},
   number={1-2},
   pages={203--240},
   issn={0944-2669},
   review={\MR{2782803}},
   doi={10.1007/s00526-010-0359-6},
}

\bib{MR3107529}{article}{
   author={Caffarelli, Luis},
   author={Valdinoci, Enrico},
   title={Regularity properties of nonlocal minimal surfaces via limiting
   arguments},
   journal={Adv. Math.},
   volume={248},
   date={2013},
   pages={843--871},
   issn={0001-8708},
   review={\MR{3107529}},
   doi={10.1016/j.aim.2013.08.007},
}

\bib{MR4000255}{article}{
   author={Cesaroni, Annalisa},
   author={Dipierro, Serena},
   author={Novaga, Matteo},
   author={Valdinoci, Enrico},
   title={Fattening and nonfattening phenomena for planar nonlocal curvature
   flows},
   journal={Math. Ann.},
   volume={375},
   date={2019},
   number={1-2},
   pages={687--736},
   issn={0025-5831},
   review={\MR{4000255}},
   doi={10.1007/s00208-018-1793-6},
}

\bib{MR3640534}{article}{
   author={Cesaroni, Annalisa},
   author={Novaga, Matteo},
   title={Volume constrained minimizers of the fractional perimeter with a
   potential energy},
   journal={Discrete Contin. Dyn. Syst. Ser. S},
   volume={10},
   date={2017},
   number={4},
   pages={715--727},
   issn={1937-1632},
   review={\MR{3640534}},
   doi={10.3934/dcdss.2017036},
}

\bib{MR3732175}{article}{
   author={Cesaroni, Annalisa},
   author={Novaga, Matteo},
   title={The isoperimetric problem for nonlocal perimeters},
   journal={Discrete Contin. Dyn. Syst. Ser. S},
   volume={11},
   date={2018},
   number={3},
   pages={425--440},
   issn={1937-1632},
   review={\MR{3732175}},
   doi={10.3934/dcdss.2018023},
}

\bib{MR3023439}{article}{
   author={Chambolle, Antonin},
   author={Morini, Massimiliano},
   author={Ponsiglione, Marcello},
   title={A nonlocal mean curvature flow and its semi-implicit time-discrete
   approximation},
   journal={SIAM J. Math. Anal.},
   volume={44},
   date={2012},
   number={6},
   pages={4048--4077},
   issn={0036-1410},
   review={\MR{3023439}},
   doi={10.1137/120863587},
}

\bib{MR3156889}{article}{
   author={Chambolle, Antonin},
   author={Morini, Massimiliano},
   author={Ponsiglione, Marcello},
   title={Minimizing movements and level set approaches to nonlocal
   variational geometric flows},
   conference={
      title={Geometric partial differential equations},
   },
   book={
      series={CRM Series},
      volume={15},
      publisher={Ed. Norm., Pisa},
   },
   date={2013},
   pages={93--104},
   review={\MR{3156889}},
   doi={10.1007/978-88-7642-473-1\_4},
}

\bib{MR3401008}{article}{
   author={Chambolle, Antonin},
   author={Morini, Massimiliano},
   author={Ponsiglione, Marcello},
   title={Nonlocal curvature flows},
   journal={Arch. Ration. Mech. Anal.},
   volume={218},
   date={2015},
   number={3},
   pages={1263--1329},
   issn={0003-9527},
   review={\MR{3401008}},
   doi={10.1007/s00205-015-0880-z},
}

\bib{MR3713894}{article}{
   author={Chambolle, Antonin},
   author={Novaga, Matteo},
   author={Ruffini, Berardo},
   title={Some results on anisotropic fractional mean curvature flows},
   journal={Interfaces Free Bound.},
   volume={19},
   date={2017},
   number={3},
   pages={393--415},
   issn={1463-9963},
   review={\MR{3713894}},
   doi={10.4171/IFB/387},
}

\bib{MR3981295}{article}{
   author={Cinti, Eleonora},
   author={Serra, Joaquim},
   author={Valdinoci, Enrico},
   title={Quantitative flatness results and $BV$-estimates for stable
   nonlocal minimal surfaces},
   journal={J. Differential Geom.},
   volume={112},
   date={2019},
   number={3},
   pages={447--504},
   issn={0022-040X},
   review={\MR{3981295}},
   doi={10.4310/jdg/1563242471},
}

\bib{MR3778164}{article}{
   author={Cinti, Eleonora},
   author={Sinestrari, Carlo},
   author={Valdinoci, Enrico},
   title={Neckpinch singularities in fractional mean curvature flows},
   journal={Proc. Amer. Math. Soc.},
   volume={146},
   date={2018},
   number={6},
   pages={2637--2646},
   issn={0002-9939},
   review={\MR{3778164}},
   doi={10.1090/proc/14002},
}

\bib{MR3836150}{article}{
   author={Ciraolo, Giulio},
   author={Figalli, Alessio},
   author={Maggi, Francesco},
   author={Novaga, Matteo},
   title={Rigidity and sharp stability estimates for hypersurfaces with
   constant and almost-constant nonlocal mean curvature},
   journal={J. Reine Angew. Math.},
   volume={741},
   date={2018},
   pages={275--294},
   issn={0075-4102},
   review={\MR{3836150}},
   doi={10.1515/crelle-2015-0088},
}

\bib{MATLUC}{article}{
   author={Cozzi, Matteo},
author={Lombardini, Luca},
title={On nonlocal minimal graphs},
      journal = {In preparation},
}

\bib{MR1942130}{article}{
   author={D\'{a}vila, J.},
   title={On an open question about functions of bounded variation},
   journal={Calc. Var. Partial Differential Equations},
   volume={15},
   date={2002},
   number={4},
   pages={519--527},
   issn={0944-2669},
   review={\MR{1942130}},
   doi={10.1007/s005260100135},
}

\bib{MR3485130}{article}{
   author={D\'{a}vila, Juan},
   author={del Pino, Manuel},
   author={Dipierro, Serena},
   author={Valdinoci, Enrico},
   title={Nonlocal Delaunay surfaces},
   journal={Nonlinear Anal.},
   volume={137},
   date={2016},
   pages={357--380},
   issn={0362-546X},
   review={\MR{3485130}},
   doi={10.1016/j.na.2015.10.009},
}

\bib{MR3412379}{article}{
   author={Di Castro, Agnese},
   author={Novaga, Matteo},
   author={Ruffini, Berardo},
   author={Valdinoci, Enrico},
   title={Nonlocal quantitative isoperimetric inequalities},
   journal={Calc. Var. Partial Differential Equations},
   volume={54},
   date={2015},
   number={3},
   pages={2421--2464},
   issn={0944-2669},
   review={\MR{3412379}},
   doi={10.1007/s00526-015-0870-x},
}

\bib{MR3007726}{article}{
   author={Dipierro, Serena},
   author={Figalli, Alessio},
   author={Palatucci, Giampiero},
   author={Valdinoci, Enrico},
   title={Asymptotics of the $s$-perimeter as $s\searrow0$},
   journal={Discrete Contin. Dyn. Syst.},
   volume={33},
   date={2013},
   number={7},
   pages={2777--2790},
   issn={1078-0947},
   review={\MR{3007726}},
   doi={10.3934/dcds.2013.33.2777},
}

\bib{MR3707346}{article}{
   author={Dipierro, Serena},
   author={Maggi, Francesco},
   author={Valdinoci, Enrico},
   title={Asymptotic expansions of the contact angle in nonlocal capillarity
   problems},
   journal={J. Nonlinear Sci.},
   volume={27},
   date={2017},
   number={5},
   pages={1531--1550},
   issn={0938-8974},
   review={\MR{3707346}},
   doi={10.1007/s00332-017-9378-1},
}

\bib{MR3516886}{article}{
   author={Dipierro, Serena},
   author={Savin, Ovidiu},
   author={Valdinoci, Enrico},
   title={Graph properties for nonlocal minimal surfaces},
   journal={Calc. Var. Partial Differential Equations},
   volume={55},
   date={2016},
   number={4},
   pages={Art. 86, 25},
   issn={0944-2669},
   review={\MR{3516886}},
   doi={10.1007/s00526-016-1020-9},
}

\bib{MR3596708}{article}{
   author={Dipierro, Serena},
   author={Savin, Ovidiu},
   author={Valdinoci, Enrico},
   title={Boundary behavior of nonlocal minimal surfaces},
   journal={J. Funct. Anal.},
   volume={272},
   date={2017},
   number={5},
   pages={1791--1851},
   issn={0022-1236},
   review={\MR{3596708}},
   doi={10.1016/j.jfa.2016.11.016},
}

\bib{2019-1}{article}{
   author={Dipierro, Serena},
   author={Savin, Ovidiu},
   author={Valdinoci, Enrico},
        title = {Nonlocal minimal graphs in the plane are generically sticky},
      journal = {arXiv e-prints},
     date = {2019},
          eid = {arXiv:1904.05393},
        pages = {arXiv:1904.05393},
archivePrefix = {arXiv},
       eprint = {1904.05393},
       adsurl = {https://ui.adsabs.harvard.edu/abs/2019arXiv190405393D},
 }

\bib{2019-2}{article}{
   author={Dipierro, Serena},
   author={Savin, Ovidiu},
   author={Valdinoci, Enrico},
        title = {Boundary properties of fractional objects: flexibility of linear equations and rigidity of minimal graphs},
      journal = {arXiv e-prints},
         date = {2019},
          eid = {arXiv:1907.01498},
        pages = {arXiv:1907.01498},
archivePrefix = {arXiv},
       eprint = {1907.01498},
       adsurl = {https://ui.adsabs.harvard.edu/abs/2019arXiv190701498D},
 }

\bib{PISA}{article}{
author={Farina, Alberto},
author={Valdinoci, Enrico},
title = {Flatness results for nonlocal minimal cones and subgraphs},
      journal = {Ann. Sc. Norm. Super. Pisa Cl. Sci. (5)},
}

\bib{MR3732178}{article}{
   author={Ferrari, Fausto},
   author={Miranda, Michele, Jr.},
   author={Pallara, Diego},
   author={Pinamonti, Andrea},
   author={Sire, Yannick},
   title={Fractional Laplacians, perimeters and heat semigroups in Carnot
   groups},
   journal={Discrete Contin. Dyn. Syst. Ser. S},
   volume={11},
   date={2018},
   number={3},
   pages={477--491},
   issn={1937-1632},
   review={\MR{3732178}},
   doi={10.3934/dcdss.2018026},
}

\bib{MR3322379}{article}{
   author={Figalli, A.},
   author={Fusco, N.},
   author={Maggi, F.},
   author={Millot, V.},
   author={Morini, M.},
   title={Isoperimetry and stability properties of balls with respect to
   nonlocal energies},
   journal={Comm. Math. Phys.},
   volume={336},
   date={2015},
   number={1},
   pages={441--507},
   issn={0010-3616},
   review={\MR{3322379}},
   doi={10.1007/s00220-014-2244-1},
}

\bib{MR3680376}{article}{
   author={Figalli, Alessio},
   author={Valdinoci, Enrico},
   title={Regularity and Bernstein-type results for nonlocal minimal
   surfaces},
   journal={J. Reine Angew. Math.},
   volume={729},
   date={2017},
   pages={263--273},
   issn={0075-4102},
   review={\MR{3680376}},
   doi={10.1515/crelle-2015-0006},
}

\bib{MR2425175}{article}{
   author={Frank, Rupert L.},
   author={Lieb, Elliott H.},
   author={Seiringer, Robert},
   title={Hardy-Lieb-Thirring inequalities for fractional Schr\"{o}dinger
   operators},
   journal={J. Amer. Math. Soc.},
   volume={21},
   date={2008},
   number={4},
   pages={925--950},
   issn={0894-0347},
   review={\MR{2425175}},
   doi={10.1090/S0894-0347-07-00582-6},
}

\bib{MR2469027}{article}{
   author={Frank, Rupert L.},
   author={Seiringer, Robert},
   title={Non-linear ground state representations and sharp Hardy
   inequalities},
   journal={J. Funct. Anal.},
   volume={255},
   date={2008},
   number={12},
   pages={3407--3430},
   issn={0022-1236},
   review={\MR{2469027}},
   doi={10.1016/j.jfa.2008.05.015},
}

\bib{MR2799577}{article}{
   author={Fusco, Nicola},
   author={Millot, Vincent},
   author={Morini, Massimiliano},
   title={A quantitative isoperimetric inequality for fractional perimeters},
   journal={J. Funct. Anal.},
   volume={261},
   date={2011},
   number={3},
   pages={697--715},
   issn={0022-1236},
   review={\MR{2799577}},
   doi={10.1016/j.jfa.2011.02.012},
}

\bib{MR2487027}{article}{
   author={Imbert, Cyril},
   title={Level set approach for fractional mean curvature flows},
   journal={Interfaces Free Bound.},
   volume={11},
   date={2009},
   number={1},
   pages={153--176},
   issn={1463-9963},
   review={\MR{2487027}},
   doi={10.4171/IFB/207},
}

\bib{Short}{article}{
author = {Julin, Vesa},
author = {La Manna, Domenico},
        title = {Short time existence of the classical solution to the
fractional mean curvature flow},
      journal = {arXiv e-prints},
     date={2019},
          eid = {arXiv:1906.10990},
        pages = {arXiv:1906.10990},
archivePrefix = {arXiv},
       eprint = {1906.10990},
       adsurl = {https://ui.adsabs.harvard.edu/abs/2019arXiv190610990J},
}

\bib{2018}{article}{
	author={Lombardini, Luca},
	title={Approximation of sets of finite fractional perimeter by smooth sets and comparison of local and global s-minimal surfaces},
	journal = {Interfaces and Free Boundaries},
    volume={20},
    date={2018},
    number={2},
    pages={261-296,},
    review={\MR{MR3827804}},
    doi={10.4171/IFB/402.}
}

\bib{MR3717439}{article}{
   author={Maggi, Francesco},
   author={Valdinoci, Enrico},
   title={Capillarity problems with nonlocal surface tension energies},
   journal={Comm. Partial Differential Equations},
   volume={42},
   date={2017},
   number={9},
   pages={1403--1446},
   issn={0360-5302},
   review={\MR{3717439}},
   doi={10.1080/03605302.2017.1358277},
}

\bib{MR3996039}{article}{
   author={Maz\'{o}n, Jos\'{e} M.},
   author={Rossi, Julio D.},
   author={Toledo, Juli\'{a}n},
   title={Nonlocal perimeter, curvature and minimal surfaces for measurable
   sets},
   journal={J. Anal. Math.},
   volume={138},
   date={2019},
   number={1},
   pages={235--279},
   issn={0021-7670},
   review={\MR{3996039}},
   doi={10.1007/s11854-019-0027-5},
}

\bib{MR1940355}{article}{
   author={Maz\cprime ya, V.},
   author={Shaposhnikova, T.},
   title={On the Bourgain, Brezis, and Mironescu theorem concerning limiting
   embeddings of fractional Sobolev spaces},
   journal={J. Funct. Anal.},
   volume={195},
   date={2002},
   number={2},
   pages={230--238},
   issn={0022-1236},
   review={\MR{1940355}},
   doi={10.1006/jfan.2002.3955},
}

\bib{PAG-C}{article}{
       author = {Pagliari, Valerio},
        title = {Halfspaces minimise nonlocal perimeter: a proof via calibrations},
      journal = {arXiv e-prints},
date={2019},
          eid = {arXiv:1905.00623},
        pages = {arXiv:1905.00623},
archivePrefix = {arXiv},
       eprint = {1905.00623},
       adsurl = {https://ui.adsabs.harvard.edu/abs/2019arXiv190500623P},
}

\bib{MR3733825}{article}{
   author={Paroni, Roberto},
   author={Podio-Guidugli, Paolo},
   author={Seguin, Brian},
   title={On the nonlocal curvatures of surfaces with or without boundary},
   journal={Commun. Pure Appl. Anal.},
   volume={17},
   date={2018},
   number={2},
   pages={709--727},
   issn={1534-0392},
   review={\MR{3733825}},
   doi={10.3934/cpaa.2018037},
}

\bib{MR3168912}{article}{
   author={Ros-Oton, Xavier},
   author={Serra, Joaquim},
   title={The Dirichlet problem for the fractional Laplacian: regularity up
   to the boundary},
   language={English, with English and French summaries},
   journal={J. Math. Pures Appl. (9)},
   volume={101},
   date={2014},
   number={3},
   pages={275--302},
   issn={0021-7824},
   review={\MR{3168912}},
   doi={10.1016/j.matpur.2013.06.003},
}

\bib{MR3090533}{article}{
   author={Savin, Ovidiu},
   author={Valdinoci, Enrico},
   title={Regularity of nonlocal minimal cones in dimension 2},
   journal={Calc. Var. Partial Differential Equations},
   volume={48},
   date={2013},
   number={1-2},
   pages={33--39},
   issn={0944-2669},
   review={\MR{3090533}},
   doi={10.1007/s00526-012-0539-7},
}

\bib{MR3951024}{article}{
   author={S\'{a}ez, Mariel},
   author={Valdinoci, Enrico},
   title={On the evolution by fractional mean curvature},
   journal={Comm. Anal. Geom.},
   volume={27},
   date={2019},
   number={1},
   pages={211--249},
   issn={1019-8385},
   review={\MR{3951024}},
   doi={10.4310/CAG.2019.v27.n1.a6},
}

\end{biblist}
\end{bibdiv}
\end{document}